\documentclass[a4paper, 10pt ]{article}
\usepackage[a4paper,left=1.75cm,right=1.75cm,top=2.5cm,bottom=2.5cm]{geometry}
\usepackage{authblk}
\usepackage[T1]{fontenc}
\usepackage[utf8]{inputenc}
\usepackage{amsmath, amssymb}
\usepackage{hyperref}
\usepackage{graphicx}
\usepackage{tabularx}
\usepackage{booktabs}
\usepackage{comment}
\DeclareUnicodeCharacter{03B2}{\ensuremath{\beta}}
\usepackage{bm}
\title{Active Subspace-Guided Free-Form Deformation with Sinkhorn Autoencoders
       for Reduced-Order Modelling of Parametrised Shape Problems }
\author[1]{Guglielmo Padula\footnote{gpadula@sissa.it}}
\author[2]{Chiara Giovannini\footnote{chgio@kth.se}}
\author[1]{Gianluigi Rozza\footnote{grozza@sissa.it}}
\author[2]{Abbas Dashtimanesh \footnote{abbasda@kth.se}}

\affil[1]{mathLab, Mathematics Area, SISSA}
\affil[2]{Department of Engineering Mechanics, KTH Royal Institute of Technology, Stockholm, Sweden}

\date{}
\begin{document}

\maketitle

\begin{abstract}

Geometrically parametrised PDEs arise in shape optimisation, where solving a PDE over many deformed domains at manageable cost is a central challenge. Free-form Deformation (FFD) parametrises shape variations through a lattice of control points, but the resulting parameter space is typically high-dimensional, with most directions barely affecting the quantity of interest (QoI). We address this with a two-stage reduction pipeline. An Active Subspace (AS)-guided structured control-point selection identifies the low-dimensional subspace of FFD weights driving the QoI. A Sinkhorn Autoencoder (SAE) is then trained on the reduced weights, learning their distribution in a compact latent space by minimising the Wasserstein distance between aggregated posterior and prior. We compare three non-intrusive Reduced-order Models (ROMs), mapping full FFD weights, AS-reduced weights, or SAE latent codes to the QoI, each using Random Forests, Gaussian Processes, and K-Nearest Neighbours. We test the methodology on a non-linear heat equation with Gaussian random Robin boundary data on a deformed Stanford Bunny, and on drag prediction for the DTC Hull bulb under FFD deformations. In both cases, the generative model captures the bulk of the QoI distribution over a comparable support, thinning the shoulder opposite the design objective, while its extremes in the design-relevant direction match the FFD sample within numerical and sampling uncertainty. ROMs learned in the SAE latent space yield smaller test errors than those on full FFD weights in most regressor-metric combinations, whereas AS-reduced weights give no comparable gain. Only the latent representation attains a positive predictive coefficient $Q^2$.
\end{abstract}

\section{Introduction}
\label{sec:introduction}

Solving partial differential equations on families of geometrically deformed domains is a recurring bottleneck in computational science and engineering. Applications include hull-form hydrodynamic optimisation~\cite{salmoiraghi_free-form_2018}, aerodynamic shape design~\cite{grey_active_2018}, patient-specific biomechanical simulation~\cite{cicci_efficient_2024} and dynamic medical imaging reconstruction on deforming anatomies~\cite{gee_digs_2026}. In all these cases, a high-fidelity solver must be called hundreds or thousands of times as the shape parameters vary, making direct simulation quickly prohibitive. Reduced-order models (ROMs) offer a principled way to address this: from a modest set of high-fidelity snapshots, they learn a surrogate that replaces expensive offline computation with cheap online evaluation~\cite{padula_brief_2024, kramer_learning_2024}.  A wide range of ROM strategies have been studied. Projection-based methods such as Proper Orthogonal Decomposition remain a standard baseline~\cite{padula_brief_2024}, while Operator Inference~\cite{kramer_learning_2024} learns structured polynomial reduced operators directly from data. Physics-informed training, which augments snapshot losses with discrete finite-element residuals, reduces the gap between data reconstruction and ROM accuracy~\cite{sibuet_discrete_2025}. Non-intrusive deep learning approaches—including fully data-driven surrogates~\cite{ivagnes_enhancing_2025, codega_machine_2026} and operator-learning methods~\cite{cicci_efficient_2024}—achieve order-of-magnitude speed-ups over classical methods.  Continuous reduced-order representations using neural fields~\cite{chang_licrom_2023} and mesh-informed neural operators~\cite{shi_mesh-informed_2025} extend these ideas to unstructured and discretisation-agnostic settings.  Parameter-dependent domains remain challenging: Buka\v{c} et al.~\cite{bukac_reduced_2024} learn domain parametrisations with deep autoencoders, while Franco et al~\cite{franco_deep_2026} use a deep orthogonal decomposition to build adaptive local bases that mitigate the Kolmogorov barrier for space-interacting parameters. Optimal transport has also become a powerful ROM tool for transported features: Khamlich et al.~\cite{khamlich_optimal_2025} employ Wasserstein-distance-based displacement interpolation to augment training data and improve tracking of moving structures, and $\beta$-variational autoencoders with transformers learn near-orthogonal latent dynamics for chaotic fluid flows~\cite{solera-rico_-variational_2024}.
 
Free-form Deformation (FFD), introduced by Sederberg and Parry~\cite{sederberg_free-form_1986}, embeds the computational domain in a tensor-product lattice of control points; moving these points smoothly deforms the enclosed geometry without re-meshing or prior knowledge of the underlying surface. This mesh-agnostic flexibility has made FFD standard in shape optimisation~\cite{salmoiraghi_free-form_2018, liu_novel_2025}, digital image correlation~\cite{chapelier_free-form_2021}, and medical image analysis~\cite{fukuda_efficient_2024}. However, practical FFD lattices have tens to hundreds of control points per spatial direction, so the parameter space quickly grows to hundreds of scalar degrees of freedom. Many of these have little effect on the Quantity of Interest (QoI), and naive regression over the full space wastes modelling capacity on directions where the response is nearly constant.  

Active Subspaces (AS) use gradients to identify a low-dimensional linear subspace of the parameter space that captures most QoI variability~\cite{grey_active_2018,constantine_active_2015}. By forming and eigendecomposing the expected outer product of the gradient, AS ranks parameter directions by average influence and projects inputs onto a few active directions. Dimensionality reduction via AS is linear, tied to the QoI gradient, and does not capture the (often nonlinear, multi-modal) distribution of shapes \cite{romor_kernelbased_2022}. 

Generative models offer a complementary perspective: rather than projecting along a gradient, they learn a latent space that captures the statistical geometry of a shape ensemble~\cite{padula_generative_2025, padula_generative_2024}. Relative to \cite{padula_generative_2024}, the principal methodological novelty consists of a reordering of the computational pipeline. In \cite{padula_generative_2024}, active subspaces are identified in the latent space learned by an autoencoder. Here, by contrast, we first compute an explicitly interpretable variant of the active subspace directly in the original input space and subsequently train the autoencoder on the resulting active-subspace representation. This modified ordering yields improved interpretability of the overall framework compared with \cite{padula_generative_2024}. Furthermore, as a novelty with respect to both \cite{padula_generative_2024} and \cite{padula_generative_2025}, we adopt Sinkhorn Autoencoders as generative models. \\
Recent deep generative methods for 3D geometry—variational autoencoders, GANs, and diffusion models on point clouds, meshes, implicit fields, and voxels~\cite{regenwetter_deep_2022, shi_deep_2023, caytuiro_3d_2025, wang_diffusion_2025}—provide a rich toolkit. Diffusion models in particular excel at 3D shape synthesis~\cite{liu_meshdiffusion_2023, xiong_octfusion_2025, potamias_shapefusion_2025, mo_efficient_2024, hu_topology-aware_2024, roessle_l3dg_2024, li_shapegen_2025, jun_shap-e_2023, siddiqui_meshgpt_2024}, and their latent spaces support controllable editing and conditional generation. Among autoencoding frameworks, the Sinkhorn Autoencoder (SAE) of Patrini et al.~\cite{patrini_sinkhorn_2020} is especially suitable. SAE minimises the $p$-Wasserstein distance between the encoder’s aggregated posterior and a latent prior by backpropagating through the Sinkhorn algorithm, operating directly on samples rather than via a reparametrization trick. This optimal-transport foundation yields strong guarantees: minimising the $p$-Wasserstein distance between the generator and the true data distribution reduces to an unconstrained min–min problem over the latent-space Wasserstein distance plus a reconstruction error~\cite{patrini_sinkhorn_2020}.

In this paper, we replace the global AS projection with an AS-guided structured control-point selection that selects FFD weights based on their marginal explanatory power from the AS matrix, keeping only those scalar components whose removal would noticeably increase the variance of the projected response, in order to increase the interpretability of the reduction.  Furthermore, we model low-dimensional selected FFD weight vectors rather than raw meshes using generative models. As a generative model, we adopt SAE, as it operates over arbitrary metric spaces and prior distributions with minimal adaptation, making it an ideal tool for learning the implicit distribution of AS-reduced FFD weights, whose geometry may differ substantially from a simple Gaussian. Once the latent representation has been learned, the SAE decoder maps any latent code back to a distribution-consistent set of FFD weights, and the associated high-fidelity solution can be obtained through the full-order solver.  This two-stage reduction---AS followed by SAE---yields three natural parameter representations for constructing non-intrusive ROMs: (i) the full FFD weight vector, (ii) the AS-projected reduced vector, and (iii) the SAE latent code. We train and compare data-driven regression models, mapping each of these representations to the target QoI, and report that on the cases considered the latent-space ROM generally attains lower test errors than the full-space one; the compactness of the latent space and its geometric regularity are advanced as possible explanations rather than as established causes. As ROMs, we adopt Random Forest (RF), Gaussian Process Regression (GPR) and K-Nearest Neighbours (KNN).

The framework is assessed on two test cases that stress complementary aspects of the pipeline. The first is a diffusion equation with a Robin boundary condition drawn from a Gaussian random function on a Stanford Bunny domain subjected to global FFD deformation. The geometry, reconstructed from range images~\cite{turk_zippered_1994}, is a complex genus-zero surface whose large FFD deformations fully explore the parameter space. The high-fidelity diffusion solver is FEniCSx~\cite{baratta_dolfinx_2023, scroggs_basix_2022}, and the QoI is the mean-squared solution, a scalar integral over the entire deformed volume. The second test case considers the drag force of the DTC Hull bulb, a canonical ship-hydrodynamics benchmark. Here, the high-fidelity solver is OpenFOAM~\cite{weller_tensorial_1998}, the FFD deformations are localised to the bulb, and the QoI depends more sharply and nonlinearly on the shape parameters than the diffusion energy. Furthermore, we restrict to the FFD weights that preserve volume. Together, these cases span diffusion-dominated and convection-dominated regimes, volumetric and surface-integral QoIs, global and localised deformations.

The remainder of this paper is organised as follows.
Section~\ref{sec:methodology} reviews the FFD parametrisation and the Active
Subspace method, describes the AS-guided structured control-point selection algorithm used to
identify the most influential FFD weights and presents the Sinkhorn Autoencoder and its application to the reduced weight distribution.
Section~\ref{sec:testcases} details the two test cases, the high-fidelity
solvers, and the experimental protocol and reports and discusses the results.
Section~\ref{sec:conclusions} concludes with remarks on limitations and future extensions.

\section{Methodology}
\label{sec:methodology}
In this section, we introduce the complete pipeline comprising Active Subspace identification, control-point selection, and the subsequent SAE–based reduction.  
We consider a reference mesh $\mathcal{M}$, defined by a set of $|\mathbf{P}|$ points $\mathbf{P}_a\in\mathbb{R}^3$, $a=1,\ldots,|\mathbf{P}|$, together with an associated graph $G$. The coordinates of each point are denoted by $(x_a,y_a,z_a)$.  
The proposed framework extends straightforwardly to spaces of arbitrary dimension; however, for clarity of exposition, we restrict attention to the three-dimensional setting and assume, without loss of generality, that $\mathbf{P}\subset[0,1]^3$.  
Our goal is to determine an optimal deformation of $\mathcal{M}$ that maximises a prescribed quantity of interest (QoI), typically evaluated via a computationally expensive numerical simulation. To mitigate the associated computational cost, we seek to construct a reduced-order model for this QoI.

To this end, we first introduce a parametrisation of the computational mesh. In this work, we adopt FFD as the deformation mapping, defined as
$$
\mathbf{FFD}_{\boldsymbol{\delta}}(\mathbf{P}_a)=
\begin{bmatrix}
\displaystyle
\sum_{i=0}^{N_x-1}\sum_{j=0}^{N_y-1}\sum_{k=0}^{N_z-1} 
\left(\frac{i}{N_x-1}+\delta_{ijk}^{x}\right) 
B_{i}^{N_x}(x_a)\,B_{j}^{N_y}(y_a)\,B_{k}^{N_z}(z_a)\\[6pt]
\displaystyle
\sum_{i=0}^{N_x-1}\sum_{j=0}^{N_y-1}\sum_{k=0}^{N_z-1} 
\left(\frac{j}{N_y-1}+\delta_{ijk}^{y}\right) 
B_{i}^{N_x}(x_a)\,B_{j}^{N_y}(y_a)\,B_{k}^{N_z}(z_a)\\[6pt]
\displaystyle
\sum_{i=0}^{N_x-1}\sum_{j=0}^{N_y-1}\sum_{k=0}^{N_z-1} 
\left(\frac{k}{N_z-1}+\delta_{ijk}^{z}\right) 
B_{i}^{N_x}(x_a)\,B_{j}^{N_y}(y_a)\,B_{k}^{N_z}(z_a)
\end{bmatrix},
$$
which is parametrised by $\delta_{ijk}^{x}$, $\delta_{ijk}^{y}$, and $\delta_{ijk}^{z}$, i.e., the Cartesian components of the displacements $\boldsymbol{\delta}_{ijk}$ of the control points from their reference positions on a uniform lattice. $B_{i}^{N_x}$ is the Bernstein polynomial of degree $N_x-1$,
$$
B_{i}^{N_{x}}(x)=\binom{N_x-1}{i}\,x^{i}(1-x)^{N_{x}-1-i},
$$
with $B_{j}^{N_y}$ and $B_{k}^{N_z}$ defined analogously. We denote by $\boldsymbol{\delta}$ the full set of control points associated with a given FFD mapping.

An important property of the FFD mapping is that it reproduces the identity transformation in the absence of deformation, that is,
\[
\mathbf{FFD}_{\boldsymbol{0}}(\mathbf{x}) = \mathbf{x}, \quad \forall\, \mathbf{x} \in [0,1]^3.
\]
We observe that the size of $\boldsymbol{\delta}$ is $N_x \cdot N_y \cdot N_z \cdot 3$, which is typically very large and may lead to difficulties in constructing a reduced-order model, owing to the curse of dimensionality~\cite{padula_generative_2024}. To mitigate this issue, parameter reduction techniques must be employed~\cite{padula_generative_2024, padula_generative_2025}. 

A classical parameter reduction method is the Active Subspaces (AS) approach. Let $\boldsymbol{\delta}^{b}$, $b = 1, \dots, M$, denote a set of control points associated with $M$ different FFD configurations, and let $f$ be the mapping that evaluates the quantity of interest on the deformed mesh induced by the FFD parameters. The AS matrix is then defined as
\[
A = \frac{1}{M} \sum_{m=1}^{M} \nabla f(\boldsymbol{\delta}^{(m)}) \, \nabla f(\boldsymbol{\delta}^{(m)})^{T}.
\]
where $\nabla f$ is flattened from a C-ordering  $(N_x,N_y,N_z,3)$. In this way, the indices of the same control point are contiguous.
The quantity $\nabla f$ is computed by kernel interpolation on $f$ and by deriving the resulting interpolator, as done in \cite{padula_generative_2024}. 

In the classical Active Subspaces (AS) framework, an eigendecomposition of the matrix \(A\) is performed, allowing us to write
\[
A = V \Lambda V^{T},
\]
where \(\Lambda\) is a diagonal matrix whose entries are the eigenvalues of \(A\), ordered in descending magnitude. The reduced parametrisation is then given by \(V_{r}^{T}\boldsymbol{\delta}\), where \(V_r\) denotes the matrix formed by the first \(r\) eigenvectors.\\
The main advantage of this technique is its computational efficiency. However, after computing the reconstructed control-point deformations \(\bar{\boldsymbol{\delta}} = \boldsymbol{V}_{r}\boldsymbol{V}_{r}^{T}\boldsymbol{\delta}\), one generally observes that \(\bar{\boldsymbol{\delta}}_{ijk}\) is nonzero even when the norm of the gradient of \(f\) with respect to \(\boldsymbol{\delta}_{ijk}\) is very small.\\
To enhance the interpretability of the reduced representation, we seek to avoid this situation.\\
The original problem can be recast as 
$$
\arg\min\limits_{V\in \mathbb{R}^{3N_xN_yN_z\times r} | V^{T}V=I} ||A-AVV^{T}||_F^{2}. 
$$
We modify it by adding a constraint
$$
\begin{array}{ll}
\arg \min \limits_{V\in\mathbb{R}^{3N_xN_yN_z\times r} } & \left\|A-AV V^T\right\|_F^{2} \\
\text { s.t. } & V \in\{0,1\}^{n \times r}, \\
&  V^{T}V=1 \\
& \sum_{i \in G_g} \sum_{k=1}^r V_{i k}=\left|G_g\right| z_g, \quad z_\alpha \in\{0,1\} , g \in 1\ldots N_xN_yN_z,
\end{array}
$$
where
$\mathcal{G}_g=\{i \in\{1, \ldots, n\}:\lfloor(i-1) / 3\rfloor=g\}$.

Thus, we force a discrete, structured selection of basis vectors in which each control-point degree of freedom is assigned in a mutually exclusive way to a single reduced component, while preserving consistency across the three displacement components belonging to the same spatial grid node. In particular, each control point can be either globally active or globally inactive.
As $V$ is composed of subsets of the identity matrix, $AVV^{T}$ has either as columns the corresponding column of $A$ or a $0$ column, we can recast the problem in a simpler way.

\begin{figure}
    \centering
    \includegraphics[width=0.75\linewidth]{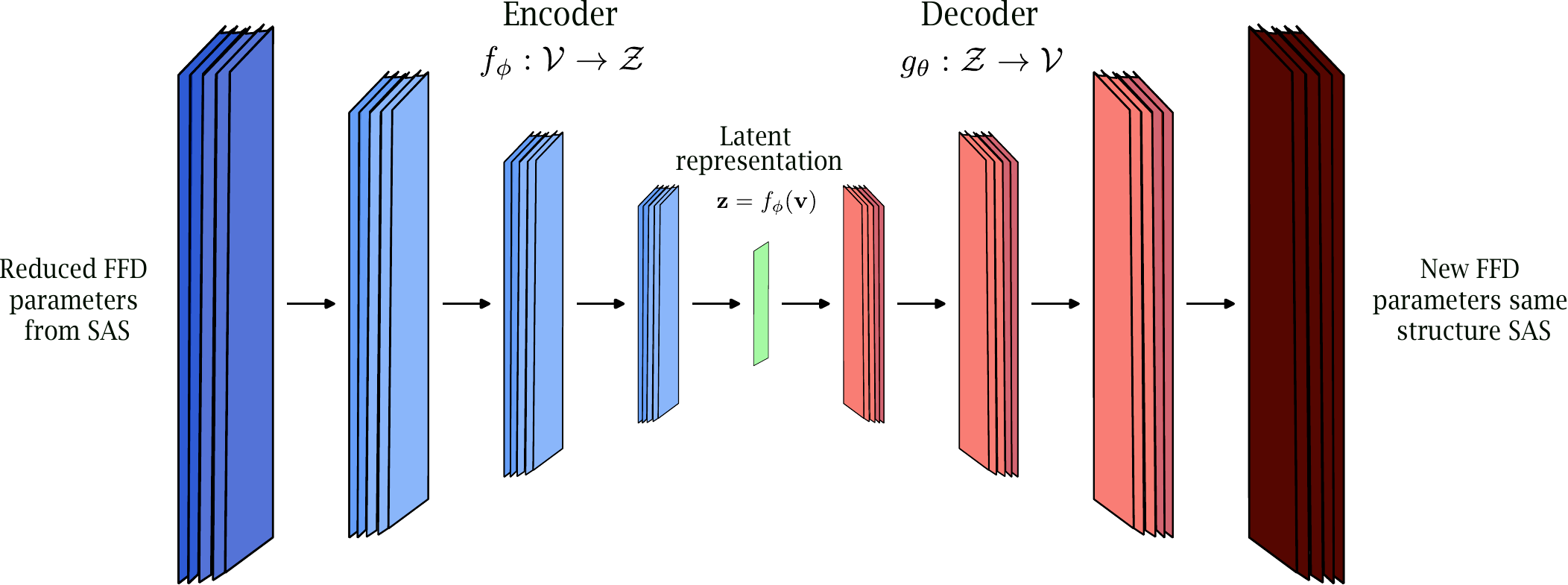}
    \caption{Architecture of the Sinkhorn Autoencoder (SAE). The deformation weights obtained through the Active Subspaces (AS) framework are provided as input to the SAE, which learns a latent representation and generates new sets of deformation weights. }
    \label{fig:SAE}
\end{figure}

To this end, let \(C_{ijk}^{x}\) the sum of the square of the column of $A$ corresponding to that index, and $C_{ijk}=C_{ijk}^{x}+C_{ijk}^{y}+C_{ijk}^{z}$.
The problem is equivalent to selecting the \(r/3\) control points corresponding to the largest values of \(C_{ijk}\). The projection matrix \(\mathbf{W}_{r}\) is then constructed by choosing as its columns the standard basis vectors associated with the selected control points. We refer to this approach as Sparse Active Subspaces (SAS), and the $0$-norm of the reconstructed FFD weights is considerably smaller than in the standard AS case. \\
A drawback of employing \(\mathbf{W}_{r}\) instead of \(\mathbf{V}_{r}\) is that a higher reduced dimensionality of the parameter space is required to achieve the same reconstruction accuracy. To address this limitation, we learn a generative model over the reduced parameterisation in order to further decrease its dimensionality. As an additional benefit, this generative model enables the synthesis of novel FFD weights, which can be advantageous in optimisation tasks.
As a generative model, we employ the Sinkhorn Autoencoder (SAE), Fig. \ref{fig:SAE}, which follows the standard encoder–decoder architecture of a conventional autoencoder~\cite{patrini_sinkhorn_2020}. Its key distinction, however, lies in the manner in which the latent space is regularised. A deterministic neural network  
\(f_\phi : \mathcal{V} \to \mathcal{Z}\)  
maps each input \(\mathbf{v} = \mathbf{W}_{r}^{T}\boldsymbol{\delta}\)
 to a latent representation  
\(\mathbf{z} = f_\phi(\mathbf{v})\). A second neural network  
\(g_\theta : \mathcal{Z} \to \mathcal{V}\)  
maps latent codes back to the reduced weight space, yielding reconstructions  
\(\hat{\mathbf{v}} = g_\theta(\mathbf{z})\).

In contrast to variational autoencoders (VAEs), which impose a factorised prior on individual latent variables, the SAE imposes a distributional constraint at the \emph{aggregate} level. Specifically, the empirical distribution of encoded training samples,
\[
  q_\phi(\mathbf{z})
  =
  \frac{1}{M}\sum_{m=1}^{M}\delta_{f_\phi(\mathbf{v}^{(m)})},
\]
is encouraged to match a standard Gaussian prior. The discrepancy between these distributions is quantified via the Sinkhorn divergence. Consequently, the overall training objective consists of a reconstruction term combined with a transport-based regularisation term:
\begin{equation}
  \mathcal{L}(\phi,\theta)
  =
  \underbrace{%
    \mathbb{E}_{\mathbf{v}\sim p_{\mathrm{data}}}
    \bigl[\,\ell\bigl(\mathbf{v},\,g_\theta(f_\phi(\mathbf{v}))\bigr)\,\bigr]
  }_{\text{reconstruction}}
  \;+\;
  \lambda\,
  \underbrace{%
    S_{\varepsilon}\!\left(q_\phi(\mathbf{z}),\,p(\mathbf{z})\right)
  }_{\text{Sinkhorn regulariser}},
  \label{eq:sae_loss}
\end{equation}
where \(\lambda > 0\) controls the trade-off between reconstruction fidelity and distributional alignment with the prior.
$S_\varepsilon$ denotes the \emph{debiased} Sinkhorn divergence, which constitutes a smoothed and computationally tractable approximation of the Wasserstein-2 (optimal transport) distance.  
Given empirical samples $\{\mathbf{z}^{(m)}\}_{m=1}^{M}$ drawn from $q_\phi$ and $\{\tilde{\mathbf{z}}_j\}_{l=1}^{L}$ drawn from $p$, $S_{\epsilon}$ is defined as
\begin{equation}
  S_\varepsilon(q,p)
  \;=\;
  \mathrm{OT}_\varepsilon(q,p)
  \;-\;
  \tfrac{1}{2}\,\mathrm{OT}_\varepsilon(q,q)
  \;-\;
  \tfrac{1}{2}\,\mathrm{OT}_\varepsilon(p,p),
  \label{eq:sinkhorn_div}
\end{equation}
where $\mathrm{OT}_\varepsilon$ denotes the entropy-regularised optimal transport cost
\begin{equation}
  \mathrm{OT}_\varepsilon(q_\phi,p)
  \;=\;
  \min_{\pi \in \Pi(q_\phi,p)}
  \sum_{m,l}
  \pi_{ij}\,c(\mathbf{z}_m, \tilde{\mathbf{z}}_l)
  \;+\;
  \varepsilon\,\mathrm{KL}(\pi \,\|\, q_\phi \otimes p),
  \label{eq:reg_ot}
\end{equation}
with $c(\cdot,\cdot)$ denoting the ground cost (typically the squared Euclidean distance) and $\varepsilon > 0$ an entropic regularisation parameter. The minimisation over couplings $\Pi(q,p)$ can be computed efficiently via the \emph{Sinkhorn--Knopp} algorithm~\cite{patrini_sinkhorn_2020}, which performs alternating row and column normalisations of the Gibbs kernel
\[
\mathbf{K}_{ml} \;=\; \exp\!\bigl(-c(\mathbf{z}_m,\tilde{\mathbf{z}}_l)/\varepsilon\bigr).
\]
The debiasing terms are chosen such that $S_\varepsilon(q,q)=0$, thereby eliminating the self-transport artefact that would otherwise introduce a systematic bias in gradient estimates~\cite{patrini_sinkhorn_2020}. 

\begin{figure}
    \centering
    \includegraphics[width=0.9\linewidth]{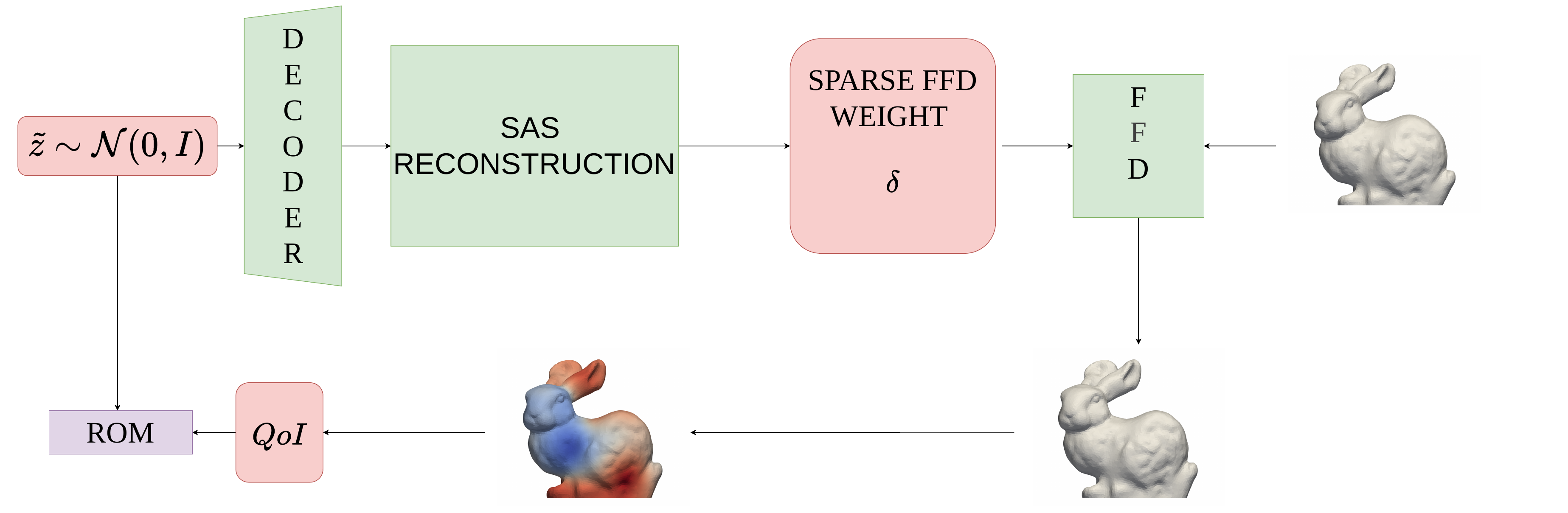}
    \caption{Detailed pipeline for the generation part and the rom assembly. Starting from a mesh $\mathcal{M}$, a free-form deformation (FFD) is applied and the quantity of interest (QoI) is evaluated. The active subspaces (AS) technique is then employed to obtain a set of deformation weights. These weights are used as input to the Sinkhorn Autoencoder (SAE), which generates a new set of weights. The corresponding QoIs are then computed and the method is finally evaluated using ROMs. }
    \label{fig:graphical_abstract}
\end{figure}

Following the training of the generative model, new sets of FFD weights are sampled, the corresponding QoIs are computed, and the mapping between the latent codes and the associated QoIs is subsequently inferred by constructing a ROM. 
A focus for the ROM assembly is provided in Fig.~\ref{fig:graphical_abstract}.

\section{Test cases and results}
\label{sec:testcases}
Our objectives are twofold: (i) to assess how closely the generative models reproduce the probability distributions of the quantities of interest (QoIs); and (ii) to compare the predictive performance of reduced-order models (ROMs) trained on the full, reduced and latent representations of the input parameters. We consider four parametrisations. Three of them describe the same dataset of FFD-generated geometries: the full FFD weight vector, the SAS-reduced coordinates obtained from it, and the SAE latent codes encoding those reduced coordinates. The fourth consists of new latent codes drawn from the generative model, which parametrise a second dataset of geometries synthesised by decoding those codes. We therefore work with two simulation datasets: one generated by FFD, on which the FFD, SAS and SAE representations are compared, and one generated by the SAE, parametrised directly in the latent space. Each dataset comprises one hundred configurations, and all the distributions reported below are empirical distributions over those hundred samples.

Two competing effects govern the accuracy of the reduced representations. On the one hand, SAS and SAE compress the parametrisation and hence cannot reproduce the original dataset exactly; on the other, the lower input dimensionality reduces the learning error of the regressor. Which effect dominates cannot be established a priori and must be assessed numerically. The latent parametrisation of the generated dataset is not affected by compression error, since the latent codes are the native parameters of the geometries in that dataset. Its distribution, however, is learned from truncated FFD coordinates and is therefore generally simpler than the original FFD distribution.

This simplification introduces an error of a different nature, which the reported metrics do not measure. A simpler input distribution yields a more regular QoI response surface, so part of the accuracy gain observed on the generated dataset may reflect an easier regression target rather than a genuinely more informative parametrisation. Errors computed on the two datasets are consequently not directly comparable, and only the FFD–SAS–SAE comparison, performed on a common set of geometries, is fully controlled. Moreover, the generative model may in principle under-represent the tails of the FFD distribution: a ROM trained on latent codes is accurate on the region covered by the generative model, but is extrapolating whenever it is queried on shapes lying outside that region, and test errors computed on generated samples cannot detect this coverage error. In the two datasets considered we find no evidence that this risk materialises in the direction relevant to design: the generated samples span supports of width comparable to the FFD ones and reach extreme values comparable to the FFD ones in the direction that an optimiser would pursue, while placing less mass on the opposite shoulder (Figs.~\ref{fig:bunny_qoi} and~\ref{fig:hull_qoi}). The caveat is nevertheless retained, since one hundred samples cannot resolve the behaviour of either distribution far into its tails. The latent-space ROM should therefore be understood as accurate with respect to the distribution it was trained on, which is the relevant one when the generative model is itself used to propose candidate shapes, as in an optimisation loop.

As reduced-order models (ROMs), we consider three commonly adopted machine-learning regression methods: Random Forests (RF), Gaussian Process Regression (GPR), and k-Nearest Neighbours (KNN). The three regressors are used with the following settings, identical for all four parametrisations and both test cases. Random Forests use $100$ trees with a maximum depth of $3$, the squared-error splitting criterion, bootstrap resampling and the remaining \texttt{scikit-learn} defaults. $k$-Nearest Neighbours uses $k=6$ with uniform weights and the Euclidean metric. Gaussian Process Regression uses a Mat\'ern kernel of smoothness $\nu=3/2$ with fixed length-scale $\ell=0.1$ and a nugget $\alpha=10^{-5}$ added to the diagonal; the hyperparameters are held at these values rather than being optimised by marginal-likelihood maximisation, so that the same kernel is used across all representations and the comparison is not confounded by differences in hyperparameter fitting.

Inputs are mapped to $[0,1]$ by a single global affine transformation, using the minimum and maximum over the whole design matrix rather than per feature, so that the relative scaling of the components of each parametrisation is preserved. The QoI is standardised to zero mean and unit variance; the reported errors are rescaled by its standard deviation and are therefore expressed in the physical units of the QoI. The normalisation constants are computed on the training set only and then applied unchanged to the test set, consistently with the treatment of the AS and SAE mappings described below.

Model performance is quantified on a single random split of each dataset into a
training set $\mathcal{T}$, comprising $80\%$ of the samples, and a disjoint test set
$\mathcal{S}$, comprising the remaining $20\%$. The split is performed before any step
of the pipeline: the active-subspace matrix, the SAS control-point selection and the
Sinkhorn Autoencoder are computed from the training samples only, and the test samples
are then projected onto the selected control points and encoded with the trained
encoder without taking part in any fitting. The reported test errors therefore assess
the generalisation performance of the complete AS--SAE--ROM workflow on unseen
geometries, and not only that of its regression stage. For the generated dataset, whose
parameters are the latent codes themselves, the split affects the regression stage
only. Since each test set contains only twenty samples, the test metrics are themselves
subject to appreciable sampling variability.

Denoting by $\hat q_m$ the prediction of the regressor fitted on $\mathcal{T}$, we
report on the test set the Mean Squared Error and the Mean Absolute Error,
$$
\text{MSE}=\frac{1}{|\mathcal{S}|}\sum_{m\in\mathcal{S}}\left(q_{m}-\hat{q}_{m}\right)^{2},
\qquad
\text{MAE}=\frac{1}{|\mathcal{S}|}\sum_{m\in\mathcal{S}}\left|q_{m}-\hat{q}_{m}\right|,
$$
together with the predictive coefficient
$$
Q^2 = 1-\frac{\sum_{m\in\mathcal{S}}(q_m-\hat q_m)^2}{\sum_{m\in\mathcal{S}}(q_m-\bar q_{\mathcal{T}})^2},
$$
where $\bar q_{\mathcal{T}}$ is the sample mean over the training set. A value
$Q^2\le0$ indicates that the regressor predicts the test set no better than the
constant predictor $\bar q_{\mathcal{T}}$.

\subsection{Nonlinear heat equation on a Stanford Bunny domain}
\label{subsec:bunny}
As a first test case, we consider a heat equation with a nonlinear diffusion coefficient subject to a Robin boundary condition on a Stanford Bunny domain~\cite{turk_zippered_1994}:
$$
\begin{cases}
\displaystyle \frac{\partial u}{\partial t}(t,\mathbf{x},\boldsymbol{\delta})
- \nabla \cdot \Bigl(\bigl(1+ u(t,\mathbf{x},\boldsymbol{\delta})^{2}\bigr)\nabla u(t,\mathbf{x},\boldsymbol{\delta})\Bigr)
= g(\mathbf{x}), 
& \mathbf{x}\in FFD_{\boldsymbol{\delta}}(\Omega),\ t\in (0,36], \\[1em]
\displaystyle u(t,\mathbf{x},\boldsymbol{\delta})
+ \bigl(1+ u(t,\mathbf{x},\boldsymbol{\delta})^{2}\bigr)\frac{\partial u}{\partial \mathbf{n}}(t,\mathbf{x},\boldsymbol{\delta}) 
= 0, 
& \mathbf{x}\in \partial FFD_{\boldsymbol{\delta}}(\Omega),\ t\in (0,36], \\[1em]
\displaystyle u(0,\mathbf{x},\boldsymbol{\delta}) = 0,
& \mathbf{x}\in FFD_{\boldsymbol{\delta}}(\Omega).
\end{cases}
$$
The solution therefore starts from rest and is driven entirely by the source term $g$.
Here, $\Omega$ denotes the undeformed Stanford Bunny domain, which is first globally rescaled to lie within $[0,1]^3$ and subsequently scaled by a factor of $0.9$ about the origin.

The initial set of $M=100$ FFD weight configurations (with resolution $N_x=N_y=N_z=5$) is sampled from a uniform distribution on $[-0.2,0.2]$, except for the control points lying on the boundary of the lattice, which are fixed to zero in order to ensure that the deformed geometries remain contained in $[0,1]^3$. Consequently, the original FFD parametrisation has $81$ effective degrees of freedom. The base of the bunny lies in the plane $z = 0$ and is therefore left undeformed. The source term $g$ is sampled once from a Gaussian random field and then held fixed for all deformed configurations.

From the active subspace matrix, we identify 16 influential control points, resulting in a reduced parametrisation with $48$ degrees of freedom, which corresponds to retaining approximately $80\%$ of the total active subspace energy.

\begin{figure}
    \centering
    \includegraphics[width=0.4\linewidth]{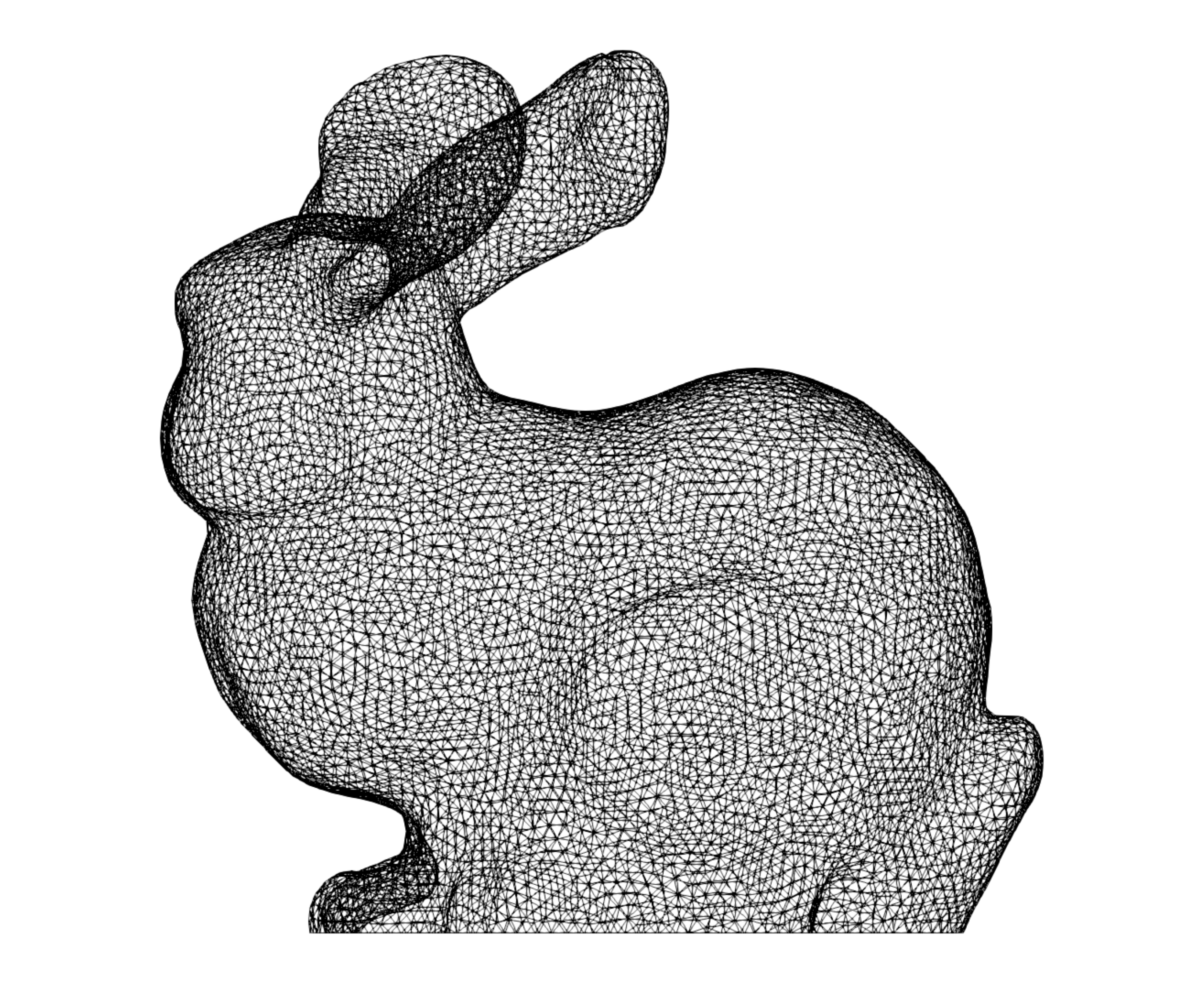} 
    \caption{The computational mesh of the undeformed Stanford Bunny.}
    \label{fig:bunny_mesh}
\end{figure}

The SAE consists of a three-layer encoder and a three-layer decoder, each with hidden layer width 500, a latent dimension of 20, ReLU activation functions, and Batch Normalisation. Training is performed using the Adam optimiser with a learning rate of $10^{-3}$ and $\epsilon = 0.01$. The bunny computational mesh is composed of $62150$ points (see Fig. \ref{fig:bunny_mesh}). The PDE is solved using the software Fenicsx~\cite{scroggs_basix_2022,baratta_dolfinx_2023} with P1 elements. As quantity of interest we consider the mean squared solution
$$
q(\boldsymbol{\delta})=\frac{1}{|\Omega_{\boldsymbol{\delta}}|}\frac{1}{T}\int_{0}^{T}\!\!\int_{\Omega_{\boldsymbol{\delta}}} u^{2}(t,\mathbf{x})\,\mathrm{d}\mathbf{x}\,\mathrm{d}t ,
\qquad \Omega_{\boldsymbol{\delta}}=FFD_{\boldsymbol{\delta}}(\Omega),
$$
which we refer to as the energy of the solution. For time stepping, we adopt a forward Euler method with $\Delta t=0.01$.

\begin{figure}
    \centering
    \includegraphics[width=0.3\linewidth]{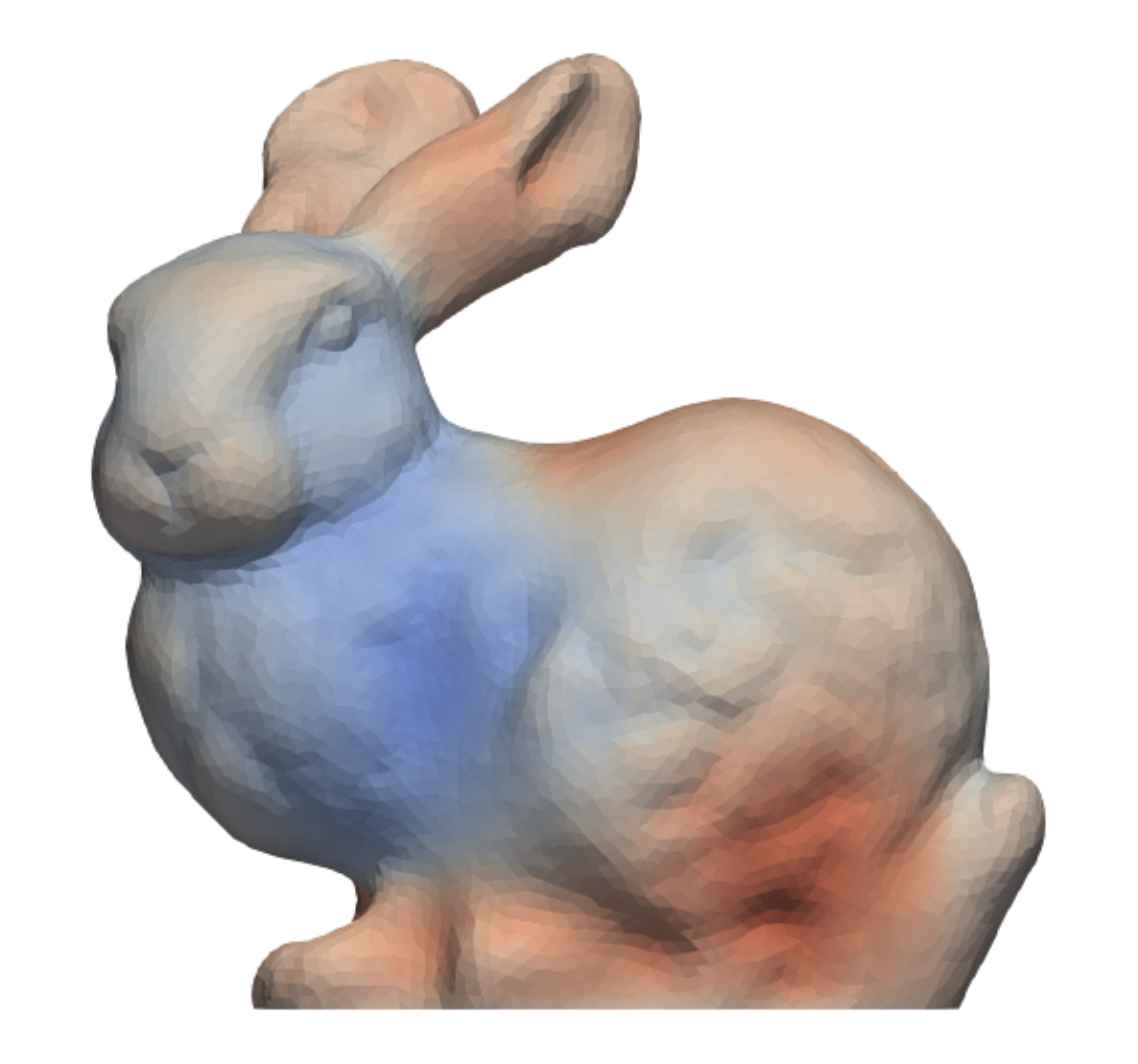}    \includegraphics[width=0.3\linewidth]{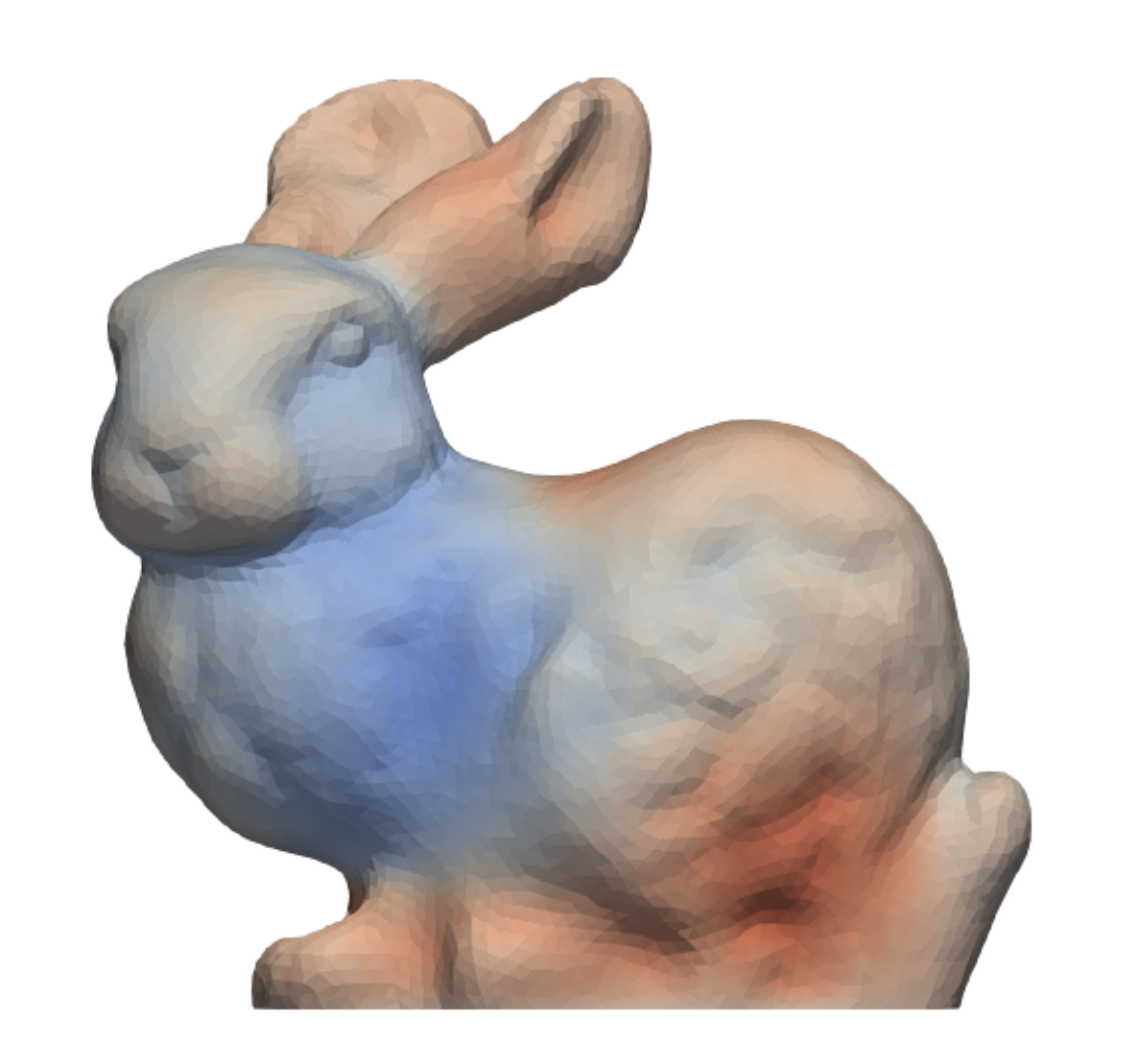}
\includegraphics[width=0.3\linewidth]{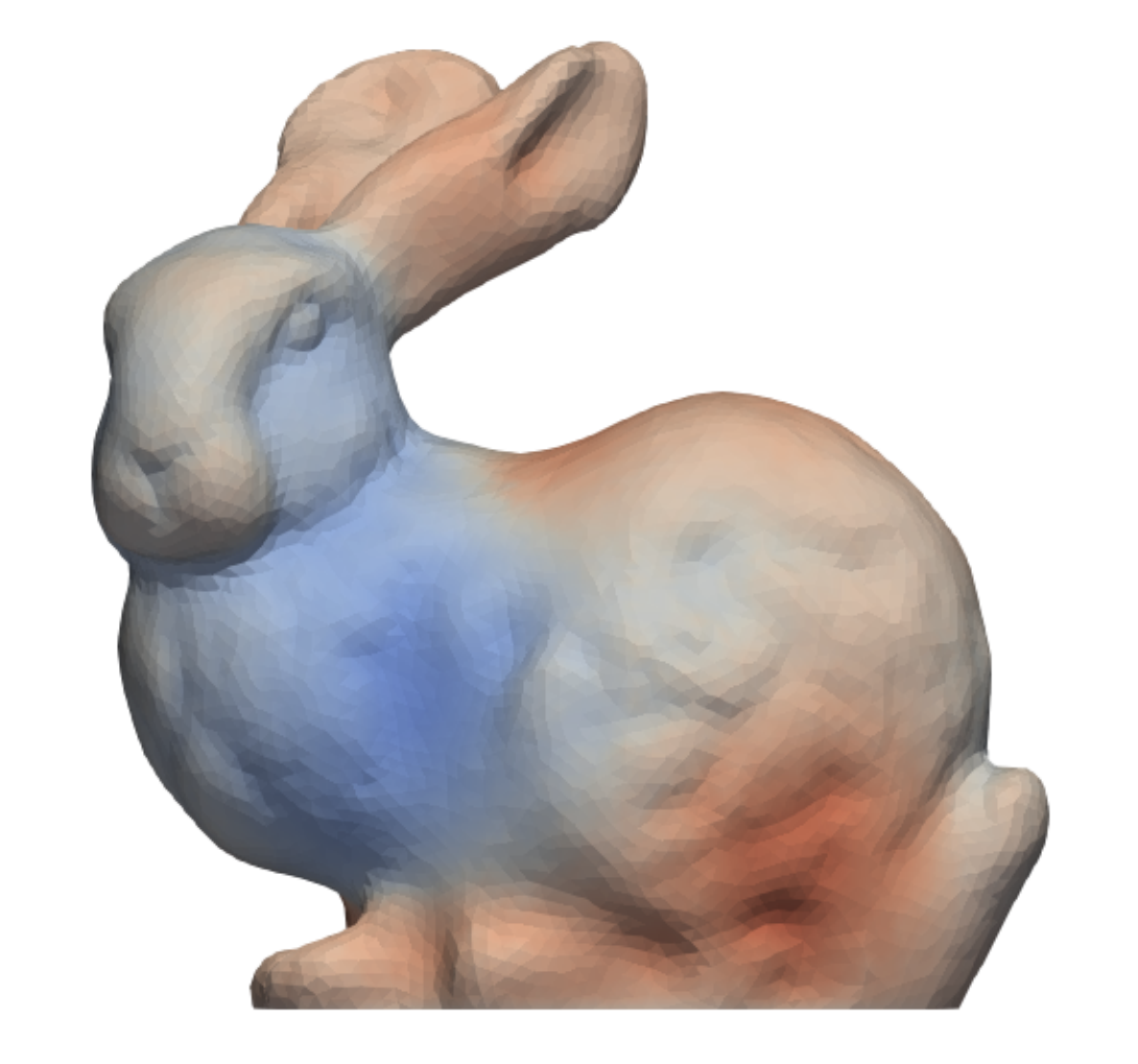}\\
\includegraphics[width=0.8\linewidth]{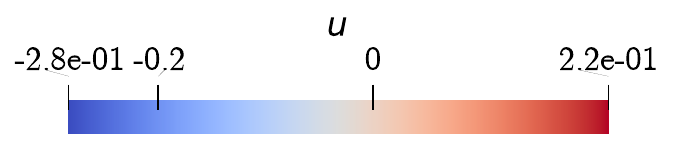}\\
    \caption{Left: solution computed on the original Stanford Bunny domain. Centre: solution computed on a domain generated via classical FFD. Right: solution computed on a domain generated using our proposed methodology. Noticeable discrepancies among the three configurations indicate a non-negligible variability in the resulting solutions.}
    \label{fig:bunny_sol}
\end{figure}

In Fig. \ref{fig:bunny_sol}, the original Stanford Bunny, an FFD-generated sample, and a sample generated by the proposed model are depicted together with their corresponding solutions, illustrating a clear variability in the resulting solutions.

\begin{figure}
    \centering
\includegraphics[width=0.75\linewidth]{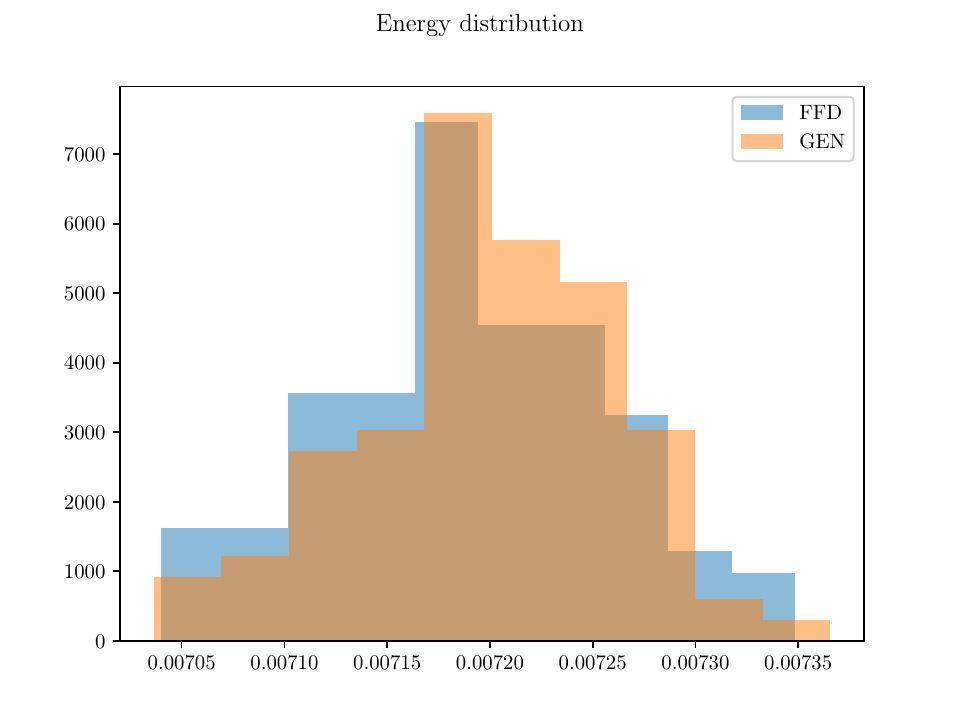}\\
\caption{Energy distribution of the FFD-generated solutions and of those produced by
the generative model, as density-normalised histograms, one hundred samples each. The
two share their location, bulk and support; the generated sample is somewhat more
concentrated, with an interquartile range of $7.9\times10^{-5}$ against
$9.5\times10^{-5}$, and carries less mass in the upper shoulder, $3\%$ above $0.0073$
against $7\%$, while its maximum is comparable to the FFD one.}
    \label{fig:bunny_qoi}
\end{figure}

Figure~\ref{fig:bunny_qoi} compares the energy distributions of the FFD-generated
solutions and of those produced by the generative model, one hundred samples each.
The two histograms agree in location, bulk and support: the mass of both is
concentrated in $[0.0071,0.0073]$, the modal bins coincide near $0.00718$, the means
differ by $1.0\times10^{-5}$, that is by $0.15$ standard deviations, and the supports,
$[0.00704,0.00735]$ for the FFD sample and $[0.00704,0.00737]$ for the generated one,
are nearly identical. The generated distribution is mildly the more concentrated of
the two, with an interquartile range of $7.9\times10^{-5}$ against $9.5\times10^{-5}$
and a standard deviation of $6.1\times10^{-5}$ against $6.7\times10^{-5}$, while the
two modal densities differ by less than $2\%$. The discrepancy is confined to the
upper shoulder, where the generated sample carries $3\%$ of its mass above $0.0073$
against $7\%$ for the FFD sample. The maxima of the two samples, $0.00737$ for the generated one and $0.00735$ for the FFD one, differ by about $0.3\%$. This difference lies within the sampling variability of the maximum of one hundred draws and has not been resolved against the discretisation error, since no mesh- or time-step-convergence study is performed; we therefore do not rank the two samples by their extremes. Were the objective to maximise the energy, the data give no indication that the generative model confines the optimiser to a range narrower than that spanned by the FFD family. The thinner upper shoulder does, however, mean that high-energy candidates are proposed less frequently per draw than under the FFD sampling distribution.

Reduced-order models yield an average computational speedup of approximately $10^{4}$
relative to the full-order simulations. This value comprises the gradient estimation
required to assemble the active subspace matrix, the SAS control-point selection, the
projection of the weights through the encoder, the fitting of the regressor and the
online evaluation, but excludes the one-off training of the Sinkhorn Autoencoder, whose
cost amounts to approximately $10^{-2}$ of the full high-fidelity dataset. The corresponding error metrics are reported in Tables~\ref{tab:mse_test_bunny}--\ref{tab:q2_test_bunny}.

For the dataset considered, the ROMs constructed in the latent representation
exhibit lower test errors than those constructed directly in the full FFD space for all
three regressors, by $23.5\%$, $32.7\%$ and $37.3\%$ in MSE and by $19.1\%$, $12.2\%$
and $13.8\%$ in MAE for RF, KNN and GPR respectively (Tables~\ref{tab:mse_test_bunny}
and~\ref{tab:mae_test_bunny}). The SAS parametrisation does not share these gains: its
MSE is higher than the full-space one by $14\%$ to $70\%$ and its MAE by $8\%$ to
$24\%$. The generative-model geometries give the lowest errors, between $27\%$ and
$51\%$ lower in MSE and between $12\%$ and $23\%$ lower in MAE than the full FFD
parametrisation, though for the reasons set out below that row is not directly
comparable with the other three.

The absolute predictive skill is, however, limited (Table~\ref{tab:q2_test_bunny}). On
the full and SAS parametrisations the test $Q^2$ is negative for all three regressors,
so none of them predicts the test energies better than the training-set mean; on the
latent representation it is positive for all three, between $0.17$ and $0.33$, and on
the generative-model geometries it reaches $0.46$ for GPR. Since neither the SAS
selection nor the SAE has seen the test geometries, these values measure the accuracy
of the whole reduction on unseen shapes rather than that of the regressor alone.

Two mechanisms are consistent with this pattern, and the present experiments do not
separate them. The first is dimensionality: at $48$ retained degrees of freedom the SAS
representation may not be compact enough for the reduction in learning error to offset
the compression error, whereas the $20$-dimensional latent space is. The second is the
regularity of the latent representation itself, since the Sinkhorn regulariser aligns
the aggregated posterior with a Gaussian prior, so that the latent coordinates are more
uniformly distributed than the FFD weights and the regression target varies more
smoothly over them. Isolating these would require holding the sampling distribution
fixed while varying the retained dimension alone, which we do not attempt here.

The lowest errors are obtained on the generative-model geometries, which carry no
compression error since the latent codes are the native parameters of those shapes. As
discussed above, part of this advantage may be attributable to the simpler input
distribution, so the four rows are not directly comparable as estimates of accuracy on
the full FFD family. They do become comparable once the intended use is shape
optimisation. In that setting the candidate geometries are themselves produced by
decoding latent codes, so the distribution on which the ROM is queried coincides with
the one on which it was trained, and no distribution mismatch arises. What matters is
then not the ability to cover the entire FFD family, but the accuracy of the
latent-to-geometry map: each latent code must decode to a well-defined admissible shape
whose QoI the ROM predicts reliably. Under this criterion the reported GEN errors are
the operationally relevant ones. This conclusion is conditional, in general, on the
generative model covering the support of the QoI distribution induced by the original
FFD family: if the decoder failed to reach a region of shape space, the corresponding
range of QoI values would be unreachable in the latent parametrisation, and the
optimiser could not recover it however accurate the ROM is on the shapes that are
reachable. Here that failure is not observed: the decoded geometries span a range of the energy
comparable to the FFD one, up to differences within the sampling and numerical
uncertainty. What is observed is a thinner upper shoulder, which lowers the rate at
which high-energy candidates are proposed without, as far as the present sample can
resolve, removing them from the reachable set; the risk remains asymmetric in
an optimisation context, because the optimum typically lies in the tail of the QoI
distribution, which is the region a generative model trained on a finite sample
resolves least well.

The dispersion of the target variable also bears on the interpretation of
Tables~\ref{tab:mse_test_bunny}--\ref{tab:q2_test_bunny}. MSE and MAE are not
normalised by the spread of the target, so a target with a smaller spread yields
smaller errors even at equal relative accuracy. This is not what distinguishes the GEN
row here: the mean squared deviations of the test targets from the training-set mean,
implied by the MSE and $Q^2$ values, are $3.31\times10^{-9}$ for the FFD geometries and
$3.20\times10^{-9}$ for the generated ones, so the advantage of the GEN row persists in
terms of $Q^2$, most markedly for GPR ($0.459$ against $0.331$ for SAE). It remains
possible that the generated configurations are simply easier to regress, for the
reasons given above, so the four rows are comparable only under the optimisation
reading just described.

\begin{table}[ht]
\caption{MSE on the test set of the ROMs considered on the energy computed on bunnies deformed using different methods. SAS denotes the set of reduced control points obtained via the active subspace matrix; SAE denotes the latent encoding via the full pipeline; GEN denotes the geometries produced by the generative model, parametrised by their latent codes.}
\label{tab:mse_test_bunny}
\centering
  \begin{tabular}{lccc}
    \hline
    Dataset & RF & KNN & GPR \\
    \hline
    FULL & $3.5996 \times 10^{-9}$ & $4.0811 \times 10^{-9}$ & $3.5330 \times 10^{-9}$ \\
    SAS  & $4.1195 \times 10^{-9}$ & $5.6741 \times 10^{-9}$ & $5.9949 \times 10^{-9}$ \\
    SAE  & $2.7551 \times 10^{-9}$ & $2.7446 \times 10^{-9}$ & $2.2155 \times 10^{-9}$ \\
    GEN  & $\bm{2.6143 \times 10^{-9}}$ & $\bm{2.4329 \times 10^{-9}}$ & $\bm{1.7296 \times 10^{-9}}$ \\
    \hline
  \end{tabular}
\end{table}

\begin{table}[ht]
\caption{MAE on the test set of the ROMs considered on the energy computed on bunnies deformed using different methods.}
\label{tab:mae_test_bunny}
\centering
  \begin{tabular}{lccc}
    \hline
    Dataset & RF & KNN & GPR \\
    \hline
    FULL & $5.1102 \times 10^{-5}$ & $4.6037 \times 10^{-5}$ & $4.7514 \times 10^{-5}$ \\
    SAS  & $5.5374 \times 10^{-5}$ & $5.6839 \times 10^{-5}$ & $5.8563 \times 10^{-5}$ \\
    SAE  & $\bm{4.1361 \times 10^{-5}}$ & $4.0429 \times 10^{-5}$ & $4.0966 \times 10^{-5}$ \\
    GEN  & $4.3253 \times 10^{-5}$ & $\bm{4.0297 \times 10^{-5}}$ & $\bm{3.6482 \times 10^{-5}}$ \\
    \hline
  \end{tabular}
\end{table}

\begin{table}[ht]
\caption{$Q^2$ on the test set of the ROMs considered on the energy computed on bunnies deformed using different methods.}
\label{tab:q2_test_bunny}
\centering
  \begin{tabular}{lccc}
    \hline
    Dataset & RF & KNN & GPR \\
    \hline
    FULL & $-0.0863$ & $-0.2316$ & $-0.0662$ \\
    SAS  & $-0.2432$ & $-0.7123$ & $-0.8091$ \\
    SAE  & $0.1686$ & $0.1717$ & $0.3314$ \\
    GEN  & $\bm{0.1826}$ & $\bm{0.2393}$ & $\bm{0.4592}$ \\
    \hline
  \end{tabular}
\end{table}

\subsection{Multiphase equation on the Duisburg test case}
\label{subsec:hull}

As a second numerical experiment, we consider a free-surface RANS simulation of the
Duisburg benchmark (DTC Hull) configuration~\cite{moctar_duisburg_2012}, as a mixed
water--air simulation.
The two phases are modelled as isothermal, immiscible and Newtonian fluids, and the
flow is described by the incompressible Reynolds-averaged Navier--Stokes equations
coupled to a volume-of-fluid transport equation for the water fraction $\alpha$. On
$\big(\Omega_{1}\cup\mathbf{FFD}_{\boldsymbol\delta}(\Omega_{2})\big)\times[0,T_{ps}]$,

\begin{align*}
& \nabla\cdot\mathbf{u} = 0, \\[1ex]
& \frac{\partial}{\partial t}\big(\rho\mathbf{u}\big)
  + \nabla\cdot\big(\rho\,\mathbf{u}\otimes\mathbf{u}\big)
  = -\nabla p_{rgh}
    - (\mathbf{g}\cdot\mathbf{x})\,\nabla\rho
    + \nabla\cdot\big(\mu\,\nabla\mathbf{u}\big)
    + \nabla\cdot R, \\[1ex]
& \frac{\partial\alpha}{\partial t}
  + \nabla\cdot\big(\mathbf{u}\,\alpha\big)
  + \nabla\cdot\big(\mathbf{u}_{r}\,\alpha(1-\alpha)\big) = 0, \\[1ex]
& \rho = \alpha\rho_W + (1-\alpha)\rho_A, \qquad
  \mu = \alpha\rho_W\nu_W + (1-\alpha)\rho_A\nu_A,
\end{align*}

where $\mathbf{u}$ is the velocity field, $\rho$ and $\mu$ the mixture density and
dynamic viscosity, $\rho_W,\rho_A$ and $\nu_W,\nu_A$ the densities and kinematic
viscosities of water and air, $R$ the Reynolds stress tensor, closed by the
$k$--$\omega$ SST model, and $\mathbf{g}$ the gravitational acceleration. The
volume fraction $\alpha$ varies between $0$ in air and $1$ in water. The mixture
kinematic viscosity, when needed, follows as $\nu=\mu/\rho$; note that the
volume-fraction weighting applies to the dynamic viscosity, so that $\nu$ is not the
linear average of $\nu_W$ and $\nu_A$.

The pressure variable is the excess pressure over the hydrostatic contribution,
\[
p_{rgh} = p - \rho\,\mathbf{g}\cdot(\mathbf{x}-\mathbf{x}_{\mathrm{ref}}),
\]
with $\mathbf{x}_{\mathrm{ref}}$ a reference point on the undisturbed free surface.
Solving for $p_{rgh}$ rather than $p$ removes the hydrostatic component from the
pressure gradient, which would otherwise have to be balanced discretely against the
buoyancy term across the density jump at the interface. This is the field reported in
Fig.~\ref{fig:hull_sol}.

The third term of the volume-fraction equation is the interface-compression flux of the
\texttt{interFoam} formulation, with $\mathbf{u}_{r}$ an artificial velocity acting
normal to the interface and scaled by the coefficient $c_{\alpha}=1$; it vanishes in
both pure phases and counteracts the numerical diffusion of the interface. Surface
tension is neglected, $\sigma=0$, as is standard for ship-resistance computations at
this scale, so no capillary term appears in the momentum equation.

A steady solution is obtained by local time stepping, the temporal derivative being
discretised with a local Euler scheme so that each cell advances at its own stable rate
and the iteration index acts as a pseudo-time rather than as physical time. Each
computation is advanced for $4000$ pseudo-time iterations, and the drag is averaged over
the last $500$; the discretisation and solver settings are given in
Appendix~\ref{app:cfd:time}. The computational mesh comprises approximately
$8.5\times10^{5}$ cells and is described in Appendix~\ref{app:cfd:mesh}.

Before introducing the deformed configurations, the full-order model is compared with
reference data for the undeformed hull. The computed drag on the half hull is
$17.0$~N, corresponding to a total resistance of $34.0$~N and, with the nominal wetted
surface $S=6.243$~m$^{2}$, to a total resistance coefficient
\[
C_{T}=\frac{R_{T}}{\tfrac{1}{2}\rho_{W}U_{0}^{2}S}=3.92\times10^{-3}.
\]
At the same Froude number, el Moctar et al.~\cite{moctar_duisburg_2012} measured
$C_{T}=3.670\times10^{-3}$ in towing-tank tests, and the potential-flow computations of
Feng et al.~\cite{feng_parametric_2021} give $3.859\times10^{-3}$. The present baseline
therefore exceeds the experimental value by $6.8\%$ and the published numerical one by
$1.6\%$. The decomposition is consistent with the expected one: the ITTC-57 correlation
line gives $C_{F}=3.05\times10^{-3}$ at $Re=9.1\times10^{6}$ which, with the form factor
$k=0.094$ reported for this hull in~\cite{feng_parametric_2021}, leaves a residuary
contribution of $5.9\times10^{-4}$, of the magnitude expected at $Fn=0.218$. The
discrepancy with the measurement is consistent with the known limitations of the
present set-up, in particular the fixed attitude and the mesh resolution inherited from
the tutorial, for which no grid-convergence study is attempted. Since the reduced-order
models are assessed on differences in drag between deformed configurations sharing the
same discretisation and the same solver settings, rather than on the absolute value of
the resistance, this level of agreement is adequate for the present purpose. We
accordingly describe the OpenFOAM computations as the full-order model of this
numerical demonstration rather than as validated high-fidelity data, and the DTC case
as an engineering demonstration of parameter-space reduction and drag-surrogate
prediction rather than as a constrained optimisation of hydrodynamically equivalent
ship designs.

The initial set of $M=100$ FFD control-point weights, associated with a lattice of
dimensions $N_x=N_y=N_z=5$ placed over the bulbous bow region (following
\cite{padula_generative_2024}), is sampled from a standard normal distribution
truncated to $[-1,1]$ and shifted by a constant offset $a=0.15$, so that each free
weight lies in $[-0.85,\,1.15]$ with mean $0.15$. The weights are then adjusted to
preserve the volume enclosed by the deformed surface, in line with
\cite{padula_generative_2024}. The offset is common to all free components, hence the
ensemble is centred on a systematically enlarged bulb rather than on the undeformed
geometry; the latter corresponds to $\boldsymbol{\delta}=\boldsymbol{0}$, which lies at
approximately the $41$st percentile of the sampling distribution and therefore remains
within the bulk of the explored configurations. Weights are expressed in units of the
lattice box, whose extent is $0.50\times0.10\times0.33$, so that before the volume
projection the mean displacement amounts to $0.075$ along $x$, $0.015$ along $y$ and
$0.050$ along $z$.

\begin{figure}
    \centering
    \includegraphics[width=0.5\linewidth]{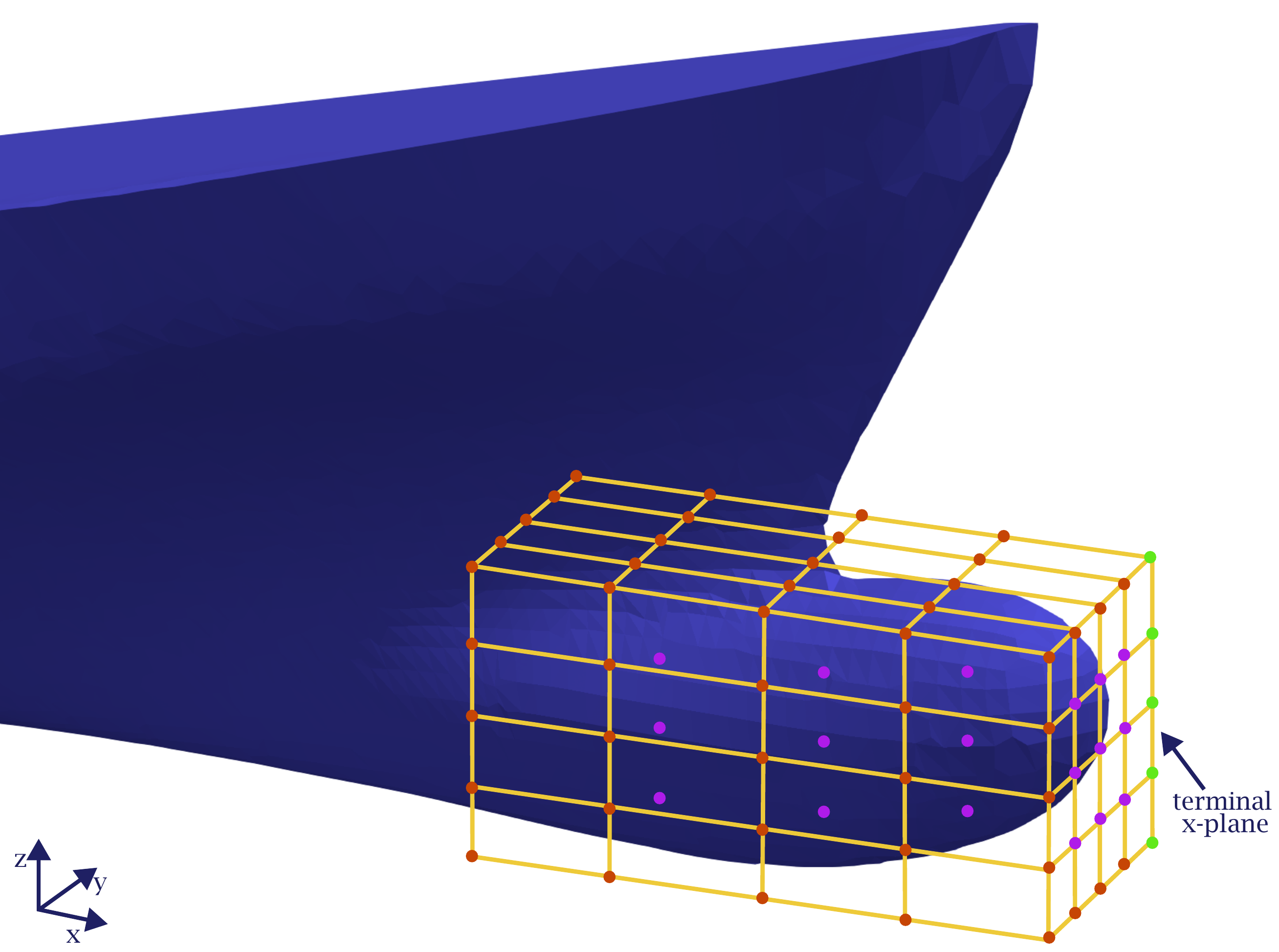}
    \caption{FFD cage in the bow region for half of the hull. Red control points denote
    the frozen outer layer, fixed on all six faces of the lattice; green points denote
    the edge released in the $x$-direction only. The control points inside the lattice
    are free in all three directions, although not all of them are displayed in the
    figure.}
    \label{fig:cage-bow}
\end{figure}

The outer layer of control points is frozen on all six faces of the lattice, so that
the deformation vanishes at the lattice boundary and the deformed bulb blends
continuously into the untouched hull, as shown in Fig.~\ref{fig:cage-bow}. This leaves
the $3\times3\times3$ interior control points free in all three directions, for a total
of $81$ degrees of freedom. In addition, the control points lying on the edge
$i=N_x-1$, $j=0$ are released in the $x$-direction only; as sampled they are assigned a
common realisation drawn as $a+10^{-2}\xi$, with $\xi$ truncated standard normal on
$[-1,1]$, so that they undergo the same mean translation as the interior points but
with a perturbation amplitude two orders of magnitude smaller. The FFD parametrisation
therefore comprises $86$ non-zero weights, of which $81$ are independent before the
volume projection; the projection distributes its correction over all free components
and relaxes the coupling of the five edge entries.

The deformation is constrained in one respect only: the control-point weights are
projected onto the constant-volume manifold, so that every configuration in the dataset
encloses the same volume as the reference bulb, to machine precision. No fairness or
curvature requirements are imposed beyond the smoothness inherent to the Bernstein
basis, and the hull is held at the fixed draught $T=0.244$~m and at fixed trim, the
computation containing no dynamic mesh. Continuity with the untouched hull is
guaranteed in position: freezing the outer layer of the lattice makes the deformation
map the identity on the lattice boundary, so the deformed bulb joins the hull without a
step; tangent continuity there is not enforced, and would require freezing two layers.
Geometries produced by the generative model are checked for validity before meshing,
and mesh quality is monitored across the dataset; the figures are reported in
Appendix~\ref{app:cfd:mesh}.

From the AS matrix, 15 control points are retained, corresponding to 45 degrees of
freedom and preserving approximately $80\%$ of the total energy; the retained weights
are again adjusted to preserve the volume, since the truncation of the discarded
components would otherwise alter it.
The SAE architecture consists of a three-layer encoder and a three-layer decoder, each
with a hidden dimension of 500, a latent dimension of 5, ReLU activation functions, and
Batch Normalisation. Training is carried out using the Adam optimiser with a learning
rate of $0.001$ and $\epsilon = 0.01$. When generating new geometries, the decoded
weights are subjected to the same volume projection.
The main quantity of interest is the drag force, measured in newtons. By convention we
define it as acting opposite to the direction of motion, so that the drag force is
taken to be positive.

\begin{figure}
    \centering
    \includegraphics[width=1.0\linewidth]{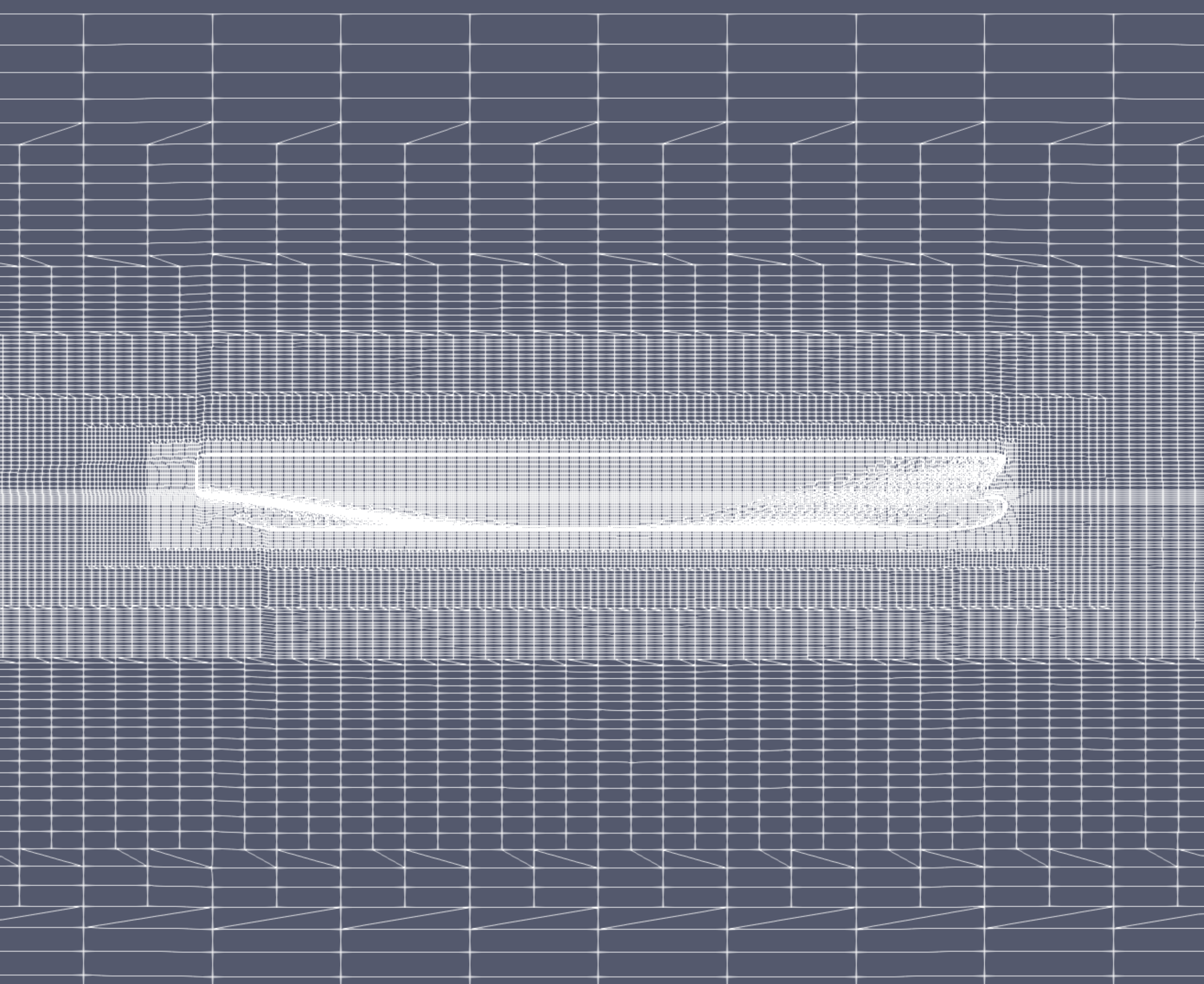}\\
    \includegraphics[width=1.0\linewidth]{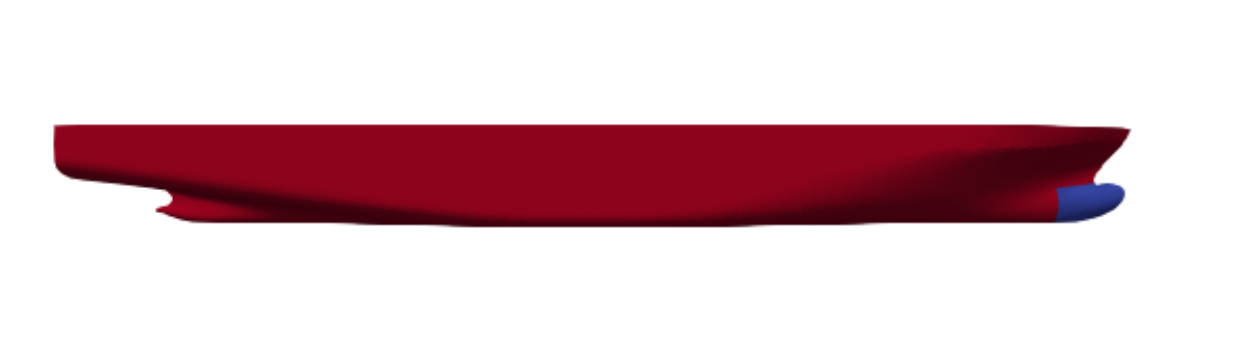}
\caption{Top: computational mesh of the undeformed DTC hull on the centreplane $y=0$.
Bottom: surface representation of the hull, with the FFD region $\Omega_{2}$ in blue
and the untouched remainder $\Omega_{1}$ in red. The grid is built in four stages: a
graded hexahedral background mesh generated with \texttt{blockMesh}, six successive
anisotropic refinements with \texttt{refineMesh} in nested boxes around the hull and
the free surface, castellation and snapping onto the triangulated hull surface with
\texttt{snappyHexMesh} (with feature-edge capture), and the addition of three
prismatic near-wall layers of expansion ratio $1.5$. The inflow is directed along
$-x$ and the calm-water plane lies at $z=0.244$~m; both are indicated. Only half of
the hull and of the domain is discretised, the centreplane being a plane of symmetry.
Full details are given in Appendix~\ref{app:cfd:mesh}.}
    \label{fig:hull_mesh}
\end{figure}

\begin{figure}
    \centering
\includegraphics[width=1\linewidth]{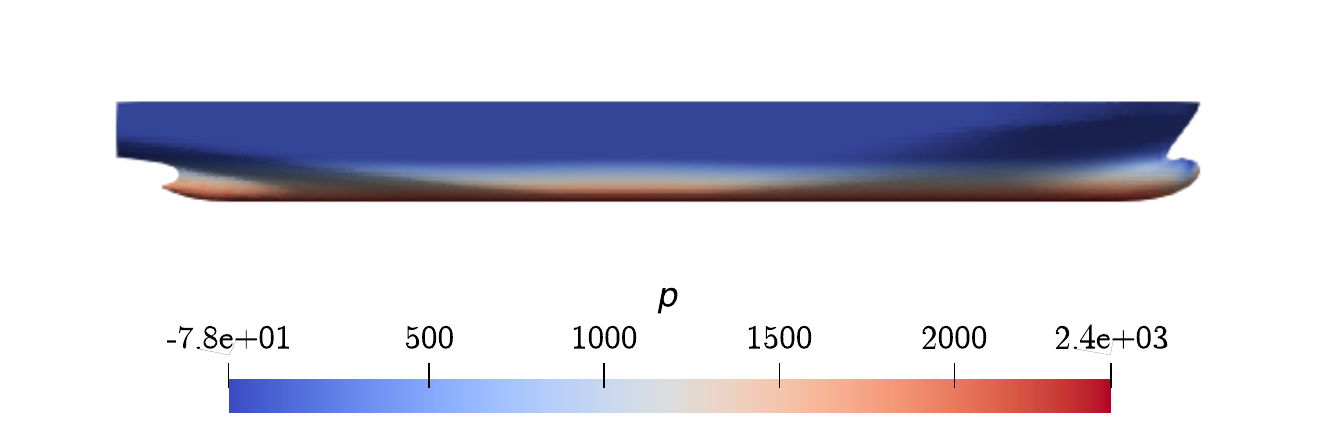}\\
\caption{Pressure field $p$ on the hull and the free surface for the baseline DTC Hull
configuration. In OpenFOAM's incompressible formulation $p$ is the kinematic pressure,
$p/\rho$, in m$^{2}$\,s$^{-2}$; the values shown are dominated by the hydrostatic
contribution $-\mathbf{g}\cdot(\mathbf{x}-\mathbf{x}_{\mathrm{ref}})$, which increases
linearly with depth below the calm-water plane $z=0.244$~m.}
    \label{fig:hull_sol_p}
\end{figure}

\begin{figure}
    \centering
    \includegraphics[width=0.3\linewidth]{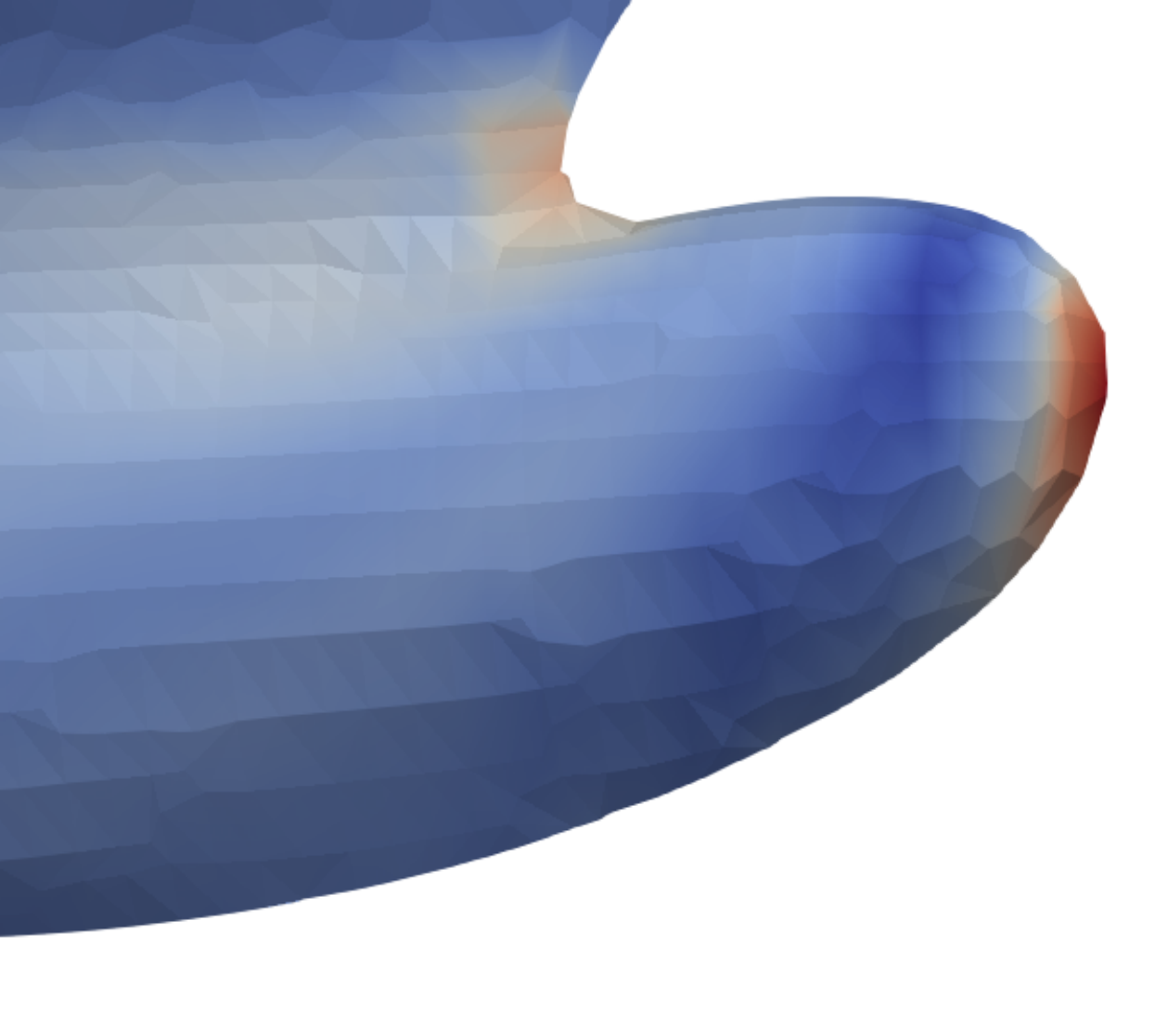}    \includegraphics[width=0.3\linewidth]{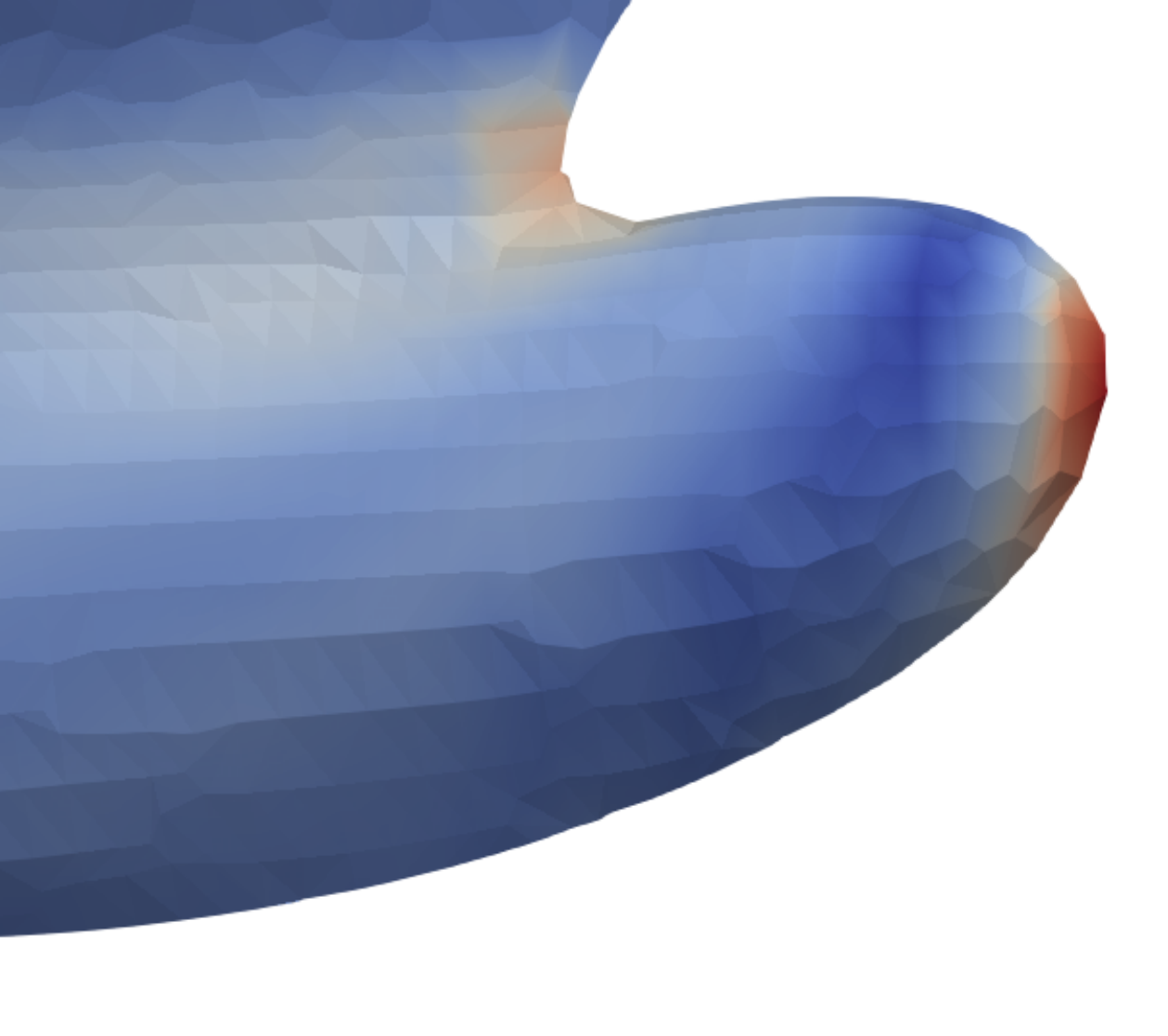}
\includegraphics[width=0.3\linewidth]{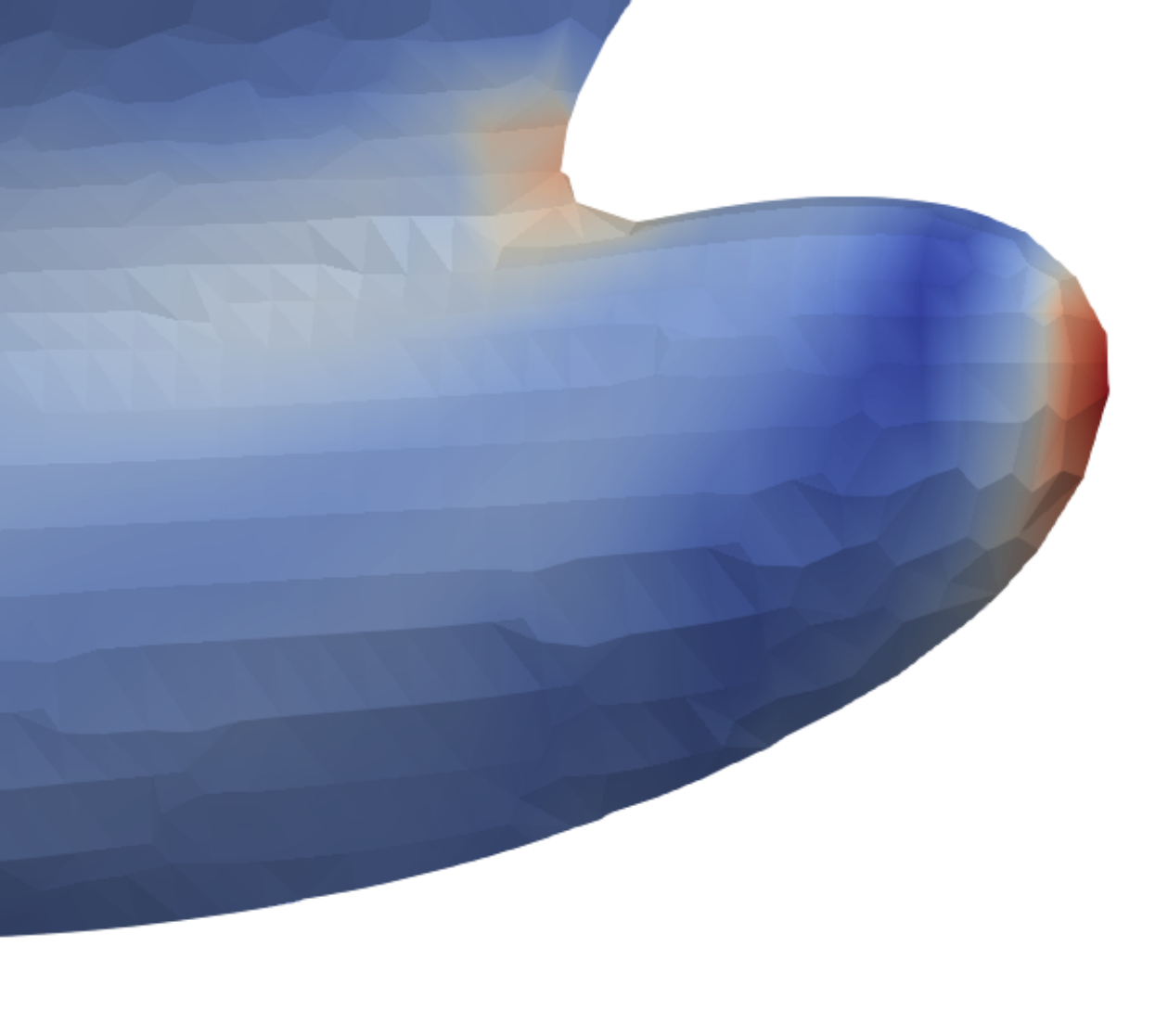}\\
\includegraphics[width=0.8\linewidth]{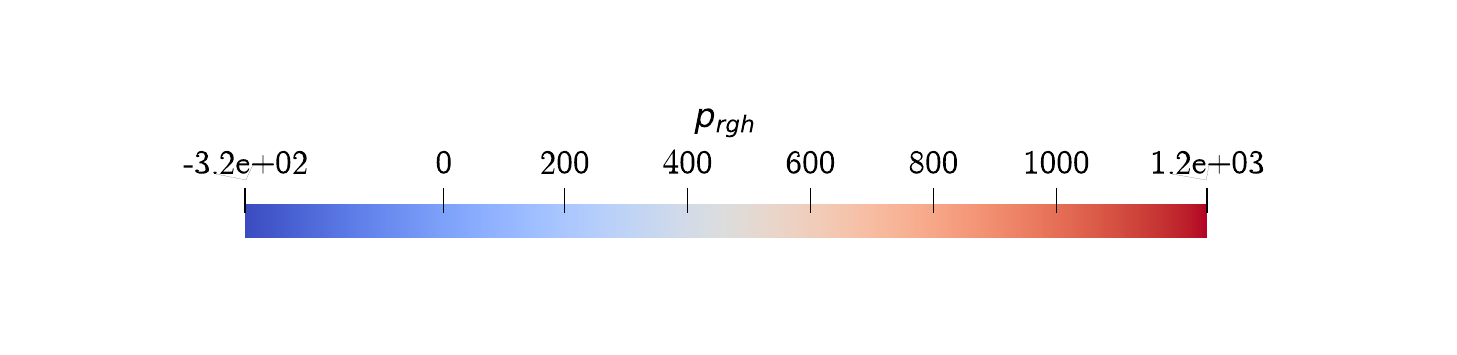}\\
    \caption{$p_{rgh}$ on the hull for the original DTC Hull domain (left), a domain
generated via classical FFD (centre), and a domain generated with the proposed
methodology (right). $p_{rgh}=p-\rho\,\mathbf{g}\cdot(\mathbf{x}-\mathbf{x}_{\mathrm{ref}})$
is the pressure with the hydrostatic part removed, in the same kinematic units; it is
the variable actually solved for, and isolates the dynamic pressure associated with
the flow, which is what differs between the three configurations. Fig.~\ref{fig:hull_sol_p}
shows $p$ for comparison. The results indicate a discernible variability among the
obtained solutions.}
    \label{fig:hull_sol}
\end{figure}

\begin{figure}
    \centering
\includegraphics[width=0.8\linewidth]{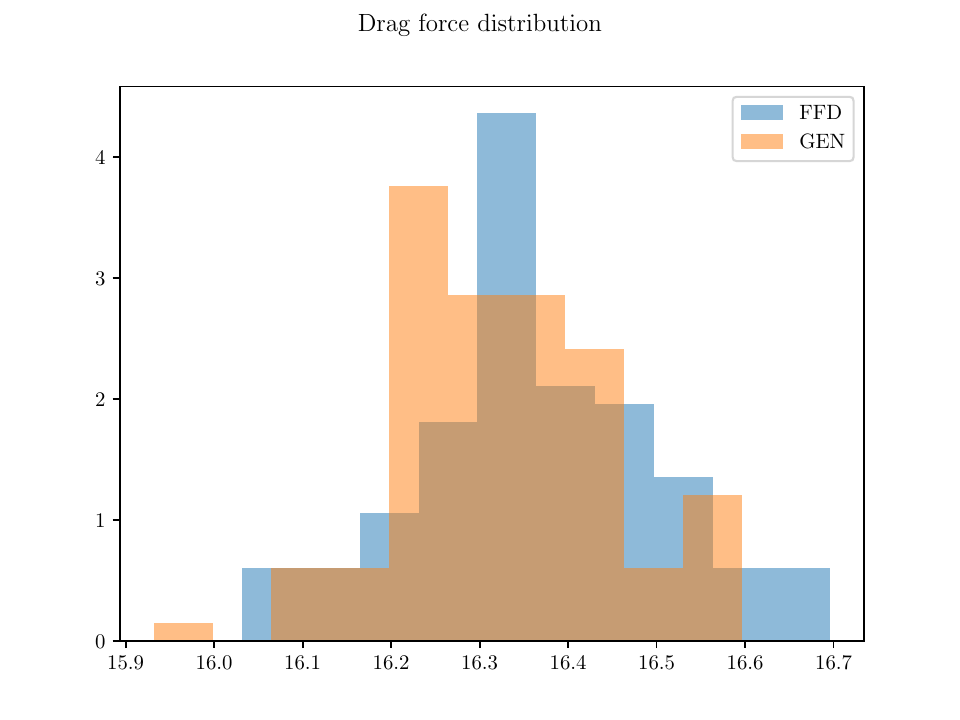}\\
\caption{Distribution of the drag force on the half hull for the FFD-generated solutions and for those produced by the generative model, as density-normalised histograms, one hundred samples
each. The two share their bulk and span supports of the same width, $0.66$~N; the
generated sample is mildly shifted toward lower drag, with a mean of $16.33$~N against
$16.36$~N, and does not exceed $16.60$~N, whereas about $6\%$ of the FFD sample lies
above that value and reaches $16.70$~N. In the direction relevant to the design problem the two minima, $15.93$~N and $16.03$~N, differ by $0.6\%$, less than the fluctuation of the force over the averaging window, and are not distinguished here.}
\label{fig:hull_qoi}
\end{figure}

Figure \ref{fig:hull_sol_p} reports the pressure field $p$ for the original DTC Hull
test case. Figure \ref{fig:hull_sol} depicts the original DTC Hull configuration, an
FFD-generated sample, and a generative-model sample, together with their corresponding
solutions for $p_{rgh}$, illustrating the variability in the resulting flow fields.
Figure \ref{fig:hull_qoi} shows the distribution of the drag force for both the
FFD-generated solutions and the generative-model solutions.

On average, the reduced-order models achieve an acceleration factor of approximately
$4\times10^{5}$ relative to the full-order simulations. This value comprises the
gradient estimation required to assemble the active subspace matrix, the SAS
control-point selection, the projection of the weights through the encoder, the fitting
of the regressor and the online evaluation, but excludes the one-off training of the
Sinkhorn Autoencoder, whose cost amounts to approximately $10^{-3}$ of the full
high-fidelity dataset. The corresponding error metrics are reported in Tables~\ref{tab:mse_test_hull}--\ref{tab:q2_test_hull}.

For the dataset considered, the ROMs constructed in the latent representation
exhibit lower test errors than those constructed directly in the full FFD space for RF
and KNN, by $18.8\%$ and $17.5\%$ in MSE and by $12.7\%$ and $13.3\%$ in MAE,
respectively (Tables~\ref{tab:mse_test_hull} and~\ref{tab:mae_test_hull}); for GPR the
MSE decreases by $3.5\%$ while the MAE increases by $15.4\%$. The SAS parametrisation
does not share these gains: its MSE lies within $3\%$ of the full-space one for all
three regressors and its MAE is higher by $2\%$ to $13\%$. The lowest errors are
attained on the generative-model geometries for all three regressors, between $58\%$
and $61\%$ lower in MSE and between $22\%$ and $32\%$ lower in MAE than the full FFD
parametrisation. The test $Q^2$ (Table~\ref{tab:q2_test_hull}) confirms the limited
absolute predictive skill observed in the diffusion case: it is negative or at most
$0.007$ on the full and SAS parametrisations, between $0.03$ and $0.13$ on the latent
one, and at most $0.23$ on the generated geometries. As in the diffusion case, the
reduced dimensionality of the latent space and its regularity under the Sinkhorn prior
are both plausible explanations for these gains, and the present comparison does not
distinguish them; that SAS, at $45$ retained degrees of freedom, does not share the gain
is consistent with the dimensionality argument. For the GEN row a third mechanism
enters, since those geometries are drawn from a distribution the decoder has learned
and are therefore not a random sample of the FFD family.
Figure~\ref{fig:hull_qoi} admits a reading that differs from the diffusion case, since
the drag force is a quantity to be minimised rather than maximised, so that the tail
relevant to the design problem is the lower one. The two distributions again share
their bulk, about $90\%$ of the mass of each lying in $[16.1,16.6]$, and the generated
sample is mildly shifted toward lower drag: its mode lies near $16.23$ against $16.33$
for the FFD sample, and its mean is $16.33$~N against $16.36$~N, a shift of $0.26$
standard deviations. The upper shoulder is the thinner of the two, as in the previous
test case: the FFD sample retains about $6\%$ of its mass above $16.60$~N and reaches
$16.70$~N, whereas no generated value exceeds $16.60$~N. Here, however, the truncation
falls in the direction irrelevant to the design problem, since the high-drag
configurations are precisely those an optimiser would discard. In the direction that matters, the two samples are comparable: their fifth percentiles
are $16.13$~N for the generated sample and $16.12$~N for the FFD one, and their minima
are $15.93$~N and $16.03$~N respectively. The latter difference, $0.10$~N or $0.6\%$,
is smaller than the $1\%$ fluctuation of the force over the averaging window
(Appendix~\ref{app:cfd:time}). It has not been assessed against the discretisation
error, for which no grid-convergence study is performed, and it concerns a single
sample out of one hundred, whose minimum is itself a highly variable statistic. We
therefore do not interpret it as the generative model proposing a configuration of
lower drag than the FFD family. The observation supports only the weaker conclusion
that the thinning of the distribution occurs on the side irrelevant to the design
problem and not on the low-drag side.

The interpretation of the error tables differs from that of the diffusion case. The
mean squared deviations of the test drag values from the training-set mean, implied by
the MSE and $Q^2$ values, are $1.96\times10^{-2}$~N$^{2}$ for the FFD geometries and
$1.00\times10^{-2}$~N$^{2}$ for the generated ones, a ratio of about one half, whereas
the variances of the full samples of Fig.~\ref{fig:hull_qoi}, $0.14^2$ and
$0.12^2$~N$^{2}$, stand in a ratio of about $0.73$; the generated test targets are thus
less dispersed than the generated sample as a whole. Once this is accounted for through
$Q^2$, the advantage of GEN over SAE largely disappears for RF ($0.141$ against $0.135$)
and KNN ($0.071$ against $0.055$), and persists only for GPR ($0.228$ against $0.034$).
The large reductions in MSE and MAE of the GEN row should therefore be attributed mainly
to the smaller spread of its test targets, which, with twenty test samples, is itself
partly a sampling effect.

A final caveat concerns the interpretation of the drag values themselves. The
deformation preserves the displaced volume exactly, so that differences in drag cannot
be attributed to changes in immersed volume, but it does not preserve the wetted
surface; since the frictional component of the resistance scales with the wetted area,
part of the variability in Fig.~\ref{fig:hull_qoi} reflects changes in surface area
rather than improvements in hydrodynamic efficiency alone. Nor is the attitude free:
redistributing the volume of the bulb moves the longitudinal centre of buoyancy, so the
deformed hulls would not in general float at the draught and trim at which they are
computed. The DTC case should accordingly be read as an engineering demonstration of
parameter-space reduction and drag-surrogate prediction, not as a constrained
optimisation over hydrodynamically equivalent ship designs; every configuration is
subject to the same constraints, so the comparison between parametrisations is
unaffected, but the ranking of configurations by drag carries no design significance.

\begin{table}
\caption{MSE (N$^2$) on the test set of the ROMs considered on the half hull drag computed on hulls deformed using different methods.}
\label{tab:mse_test_hull}
\centering
  \begin{tabular}{lccc}
    \hline
    Dataset & RF & KNN & GPR \\
    \hline
    FULL & $2.0844 \times 10^{-2}$ & $2.2426 \times 10^{-2}$ & $1.9606 \times 10^{-2}$ \\
    SAS  & $2.0224 \times 10^{-2}$ & $2.2303 \times 10^{-2}$ & $1.9440 \times 10^{-2}$ \\
    SAE  & $1.6933 \times 10^{-2}$ & $1.8503 \times 10^{-2}$ & $1.8913 \times 10^{-2}$ \\
    GEN  & $\bm{8.6084 \times 10^{-3}}$ & $\bm{9.3050 \times 10^{-3}}$ & $\bm{7.7404 \times 10^{-3}}$ \\
    \hline
  \end{tabular}
\end{table}

\begin{table}
\caption{MAE (N) on the test set of the ROMs considered on the half hull drag computed on hulls deformed using different methods.}
\label{tab:mae_test_hull}
\centering
  \begin{tabular}{lccc}
    \hline
    Dataset & RF & KNN & GPR \\
    \hline
    FULL & $1.0429 \times 10^{-1}$ & $1.1834 \times 10^{-1}$ & $9.8130 \times 10^{-2}$ \\
    SAS  & $1.0651 \times 10^{-1}$ & $1.2296 \times 10^{-1}$ & $1.1076 \times 10^{-1}$ \\
    SAE  & $9.1013 \times 10^{-2}$ & $1.0263 \times 10^{-1}$ & $1.1327 \times 10^{-1}$ \\
    GEN  & $\bm{7.1032 \times 10^{-2}}$ & $\bm{8.1610 \times 10^{-2}}$ & $\bm{7.6020 \times 10^{-2}}$ \\
    \hline
  \end{tabular}
\end{table}

\begin{table}
\caption{$Q^2$ on the test set of the ROMs considered on the half hull drag computed on hulls deformed using different methods.}
\label{tab:q2_test_hull}
\centering
  \begin{tabular}{lccc}
    \hline
    Dataset & RF & KNN & GPR \\
    \hline
    FULL & $-0.0651$ & $-0.1460$ & $-0.0019$ \\
    SAS  & $-0.2432$ & $-0.1397$ & $0.0066$ \\
    SAE  & $0.1347$ & $0.0545$ & $0.0335$ \\
    GEN  & $\bm{0.1408}$ & $\bm{0.0713}$ & $\bm{0.2275}$ \\
    \hline
  \end{tabular}
\end{table}
\section{Conclusion}
\label{sec:conclusions}
This work introduced a framework for parameter-space reduction of parametrised shape problems, combining Active Subspaces, control-point selection, and Sinkhorn Autoencoder, with applications to ROM. The proposed methodology first identifies the most influential Free-form Deformation parameters, then learns a compact latent representation of the resulting reduced parameter space, and finally exploits this representation to construct non-intrusive ROMs.  

Numerical experiments on a nonlinear diffusion problem defined on a deformed Stanford
Bunny geometry and on the DTC Hull bulb benchmark indicate that the Sinkhorn
Autoencoder reproduces the bulk of the empirical distribution of the quantities of
interest over a support of comparable width, thinning the shoulder that lies away from the design objective without a detectable reduction of the attainable range in the direction relevant to design, and that for the datasets considered the ROMs built on the learned latent representation attain lower test errors than those built on the original Free-form Deformation parameters in most regressor--metric combinations, and are the only ones with a positive test $Q^2$ for all regressors, whereas the variables obtained through Active Subspaces yield no comparable gain. The absolute predictive skill nevertheless remains modest, with test $Q^2$ not exceeding $0.33$ on the latent representation and $0.46$ overall, so the reported gains are relative improvements over regressors of limited accuracy rather than evidence of a reliable surrogate. The compactness of the latent space and its
regularity under the Sinkhorn prior are plausible explanations for this behaviour,
which the present experiments do not separate. Substantial computational speedups are
preserved throughout.

These findings suggest that the combination of sensitivity-based parameter selection with optimal-transport-based generative modelling constitutes a promising strategy for mitigating the complexity of high-dimensional shape-parameterisation problems. Future research will focus on larger training sets to raise the absolute predictive accuracy, on extending the framework to multi-objective optimisation, higher-dimensional parameter spaces, direct reduced-order modelling of field solutions rather than solely scalar quantities of interest, and on experiments designed to separate the effect of the retained dimension, together with grid-convergence and repeated-sampling studies to establish whether the generative model can extend the attainable range of the quantity of interest beyond that of the Free-form Deformation family.

\bibliographystyle{unsrt}
\bibliography{biblio}

@inproceedings{patrini_sinkhorn_2020,
	title = {Sinkhorn {AutoEncoders}},
	url = {https://proceedings.mlr.press/v115/patrini20a.html},
	language = {en},
	urldate = {2026-06-03},
	booktitle = {Proceedings of {The} 35th {Uncertainty} in {Artificial} {Intelligence} {Conference}},
	publisher = {PMLR},
	author = {Patrini, Giorgio and Berg, Rianne van den and Forré, Patrick and Carioni, Marcello and Bhargav, Samarth and Welling, Max and Genewein, Tim and Nielsen, Frank},
	month = aug,
	year = {2020},
	pages = {733--743},
}

@incollection{gee_digs_2026,
	address = {Cham},
	title = {{DIGS}: {Dynamic} {CBCT} {Reconstruction} {Using} {Deformation}-{Informed} {4D} {Gaussian} {Splatting} and a {Low}-{Rank} {Free}-{Form} {Deformation} {Model}},
	volume = {15963},
	isbn = {9783032049643 9783032049650},
	shorttitle = {{DIGS}},
	url = {https://link.springer.com/10.1007/978-3-032-04965-0_13},
	language = {en},
	urldate = {2026-06-03},
	booktitle = {Medical {Image} {Computing} and {Computer} {Assisted} {Intervention} – {MICCAI} 2025},
	publisher = {Springer Nature Switzerland},
	author = {Huang, Yuliang and Singh, Imraj and Joyce, Thomas and Thielemans, Kris and McClelland, Jamie R.},
	editor = {Gee, James C. and Alexander, Daniel C. and Hong, Jaesung and Iglesias, Juan Eugenio and Sudre, Carole H. and Venkataraman, Archana and Golland, Polina and Kim, Jong Hyo and Park, Jinah},
	year = {2026},
	doi = {10.1007/978-3-032-04965-0_13},
	pages = {131--141},
}

@article{codega_machine_2026,
	title = {Machine {Learning}-{Based} {Quadratic} {Closures} for {Non}-{Intrusive} {Reduced} {Order} {Models}},
	volume = {48},
	issn = {1064-8275, 1095-7197},
	url = {https://epubs.siam.org/doi/10.1137/25M1766759},
	doi = {10.1137/25M1766759},
	language = {en},
	number = {3},
	urldate = {2026-06-03},
	journal = {SIAM Journal on Scientific Computing},
	author = {Codega, Gabriele and Ivagnes, Anna and Demo, Nicola and Rozza, Gianluigi},
	month = jun,
	year = {2026},
	pages = {C505--C525},
}

@article{khamlich_optimal_2025,
	title = {Optimal transport-based displacement interpolation with data augmentation for reduced order modeling of nonlinear dynamical systems},
	volume = {531},
	issn = {00219991},
	url = {https://linkinghub.elsevier.com/retrieve/pii/S0021999125002219},
	doi = {10.1016/j.jcp.2025.113938},
	language = {en},
	urldate = {2026-06-03},
	journal = {Journal of Computational Physics},
	author = {Khamlich, Moaad and Pichi, Federico and Girfoglio, Michele and Quaini, Annalisa and Rozza, Gianluigi},
	month = jun,
	year = {2025},
	pages = {113938},
}

@article{ivagnes_enhancing_2025,
	title = {Enhancing non-intrusive reduced-order models with space-dependent aggregation methods},
	volume = {236},
	issn = {0001-5970, 1619-6937},
	url = {https://link.springer.com/10.1007/s00707-024-04007-9},
	doi = {10.1007/s00707-024-04007-9},
	language = {en},
	number = {9},
	urldate = {2026-06-03},
	journal = {Acta Mechanica},
	author = {Ivagnes, Anna and Tonicello, Niccolò and Cinnella, Paola and Rozza, Gianluigi},
	month = sep,
	year = {2025},
	pages = {5875--5904},
}

@article{kramer_learning_2024,
	title = {Learning {Nonlinear} {Reduced} {Models} from {Data} with {Operator} {Inference}},
	volume = {56},
	copyright = {http://creativecommons.org/licenses/by/4.0/},
	issn = {0066-4189, 1545-4479},
	url = {https://www.annualreviews.org/doi/10.1146/annurev-fluid-121021-025220},
	doi = {10.1146/annurev-fluid-121021-025220},
	language = {en},
	number = {1},
	urldate = {2026-06-03},
	journal = {Annual Review of Fluid Mechanics},
	author = {Kramer, Boris and Peherstorfer, Benjamin and Willcox, Karen E.},
	month = jan,
	year = {2024},
	pages = {521--548},
}

@article{solera-rico_-variational_2024,
	title = {β-{Variational} autoencoders and transformers for reduced-order modelling of fluid flows},
	volume = {15},
	issn = {2041-1723},
	url = {https://www.nature.com/articles/s41467-024-45578-4},
	doi = {10.1038/s41467-024-45578-4},
	language = {en},
	number = {1},
	urldate = {2026-06-03},
	journal = {Nature Communications},
	author = {Solera-Rico, Alberto and Sanmiguel Vila, Carlos and Gómez-López, Miguel and Wang, Yuning and Almashjary, Abdulrahman and Dawson, Scott T. M. and Vinuesa, Ricardo},
	month = feb,
	year = {2024},
	pages = {1361},
}

@article{cicci_efficient_2024,
	title = {Efficient approximation of cardiac mechanics through reduced‐order modeling with deep learning‐based operator approximation},
	volume = {40},
	issn = {2040-7939, 2040-7947},
	url = {https://onlinelibrary.wiley.com/doi/10.1002/cnm.3783},
	doi = {10.1002/cnm.3783},
	language = {en},
	number = {1},
	urldate = {2026-06-03},
	journal = {International Journal for Numerical Methods in Biomedical Engineering},
	author = {Cicci, Ludovica and Fresca, Stefania and Manzoni, Andrea and Quarteroni, Alfio},
	month = jan,
	year = {2024},
	pages = {e3783},
}

@article{franco_deep_2026,
	title = {Deep orthogonal decomposition: a continuously adaptive neural network approach to model order reduction of parametrized partial differential equations},
	volume = {52},
	issn = {1019-7168, 1572-9044},
	shorttitle = {Deep orthogonal decomposition},
	url = {https://link.springer.com/10.1007/s10444-026-10295-7},
	doi = {10.1007/s10444-026-10295-7},
	language = {en},
	number = {3},
	urldate = {2026-06-03},
	journal = {Advances in Computational Mathematics},
	author = {Franco, Nicola Rares and Manzoni, Andrea and Zunino, Paolo and Hesthaven, Jan S.},
	month = jun,
	year = {2026},
	pages = {31},
}

@misc{sibuet_discrete_2025,
	title = {A discrete physics-informed training for projection-based reduced order models with neural networks},
	url = {https://arxiv.org/abs/2504.13875v2},
	language = {en},
	urldate = {2026-06-03},
	journal = {arXiv.org},
	author = {Sibuet, N. and de Parga, S. Ares and Bravo, J. R. and Rossi, R.},
	month = mar,
	year = {2025},
	doi = {10.3390/axioms14050385},
}

@misc{bukac_reduced_2024,
	title = {Reduced {Order} {Modeling} of {Partial} {Differential} {Equations} on {Parameter}-{Dependent} {Domains} {Using} {Deep} {Neural} {Networks}},
	url = {https://arxiv.org/abs/2407.17171v2},
	language = {en},
	urldate = {2026-06-03},
	journal = {arXiv.org},
	author = {Bukač, Martina and Manojlović, Iva and Muha, Boris and Vlah, Domagoj},
	month = jul,
	year = {2024},
}

@inproceedings{chang_licrom_2023,
	address = {Sydney NSW Australia},
	title = {{LiCROM}: {Linear}-{Subspace} {Continuous} {Reduced} {Order} {Modeling} with {Neural} {Fields}},
	isbn = {9798400703157},
	shorttitle = {{LiCROM}},
	url = {https://dl.acm.org/doi/10.1145/3610548.3618158},
	doi = {10.1145/3610548.3618158},
	language = {en},
	urldate = {2026-06-03},
	booktitle = {{SIGGRAPH} {Asia} 2023 {Conference} {Papers}},
	publisher = {ACM},
	author = {Chang, Yue and Chen, Peter Yichen and Wang, Zhecheng and Chiaramonte, Maurizio M. and Carlberg, Kevin and Grinspun, Eitan},
	month = dec,
	year = {2023},
	pages = {1--12},
}

@inproceedings{li_shapegen_2025,
	address = {Hong Kong Hong Kong},
	title = {{ShapeGen}: {Towards} {High}-{Quality} {3D} {Shape} {Synthesis}},
	isbn = {9798400721373},
	shorttitle = {{ShapeGen}},
	url = {https://dl.acm.org/doi/10.1145/3757377.3763812},
	doi = {10.1145/3757377.3763812},
	language = {en},
	urldate = {2026-06-03},
	booktitle = {Proceedings of the {SIGGRAPH} {Asia} 2025 {Conference} {Papers}},
	publisher = {ACM},
	author = {Li, Yangguang and He, Xianglong and Zou, Zi-Xin and Liu, Zexiang and Ouyang, Wanli and Liang, Ding and Cao, Yan-Pei},
	month = dec,
	year = {2025},
	pages = {1--12},
}

@inproceedings{roessle_l3dg_2024,
	address = {Tokyo Japan},
	title = {{L3DG}: {Latent} {3D} {Gaussian} {Diffusion}},
	isbn = {9798400711312},
	shorttitle = {{L3DG}},
	url = {https://dl.acm.org/doi/10.1145/3680528.3687699},
	doi = {10.1145/3680528.3687699},
	language = {en},
	urldate = {2026-06-03},
	booktitle = {{SIGGRAPH} {Asia} 2024 {Conference} {Papers}},
	publisher = {ACM},
	author = {Roessle, Barbara and Müller, Norman and Porzi, Lorenzo and Rota Bulò, Samuel and Kontschieder, Peter and Dai, Angela and Nießner, Matthias},
	month = dec,
	year = {2024},
	pages = {1--11},
}

@misc{hu_topology-aware_2024,
	title = {Topology-{Aware} {Latent} {Diffusion} for {3D} {Shape} {Generation}},
	url = {https://arxiv.org/abs/2401.17603v1},
	language = {en},
	urldate = {2026-06-03},
	journal = {arXiv.org},
	author = {Hu, Jiangbei and Fei, Ben and Xu, Baixin and Hou, Fei and Yang, Weidong and Wang, Shengfa and Lei, Na and Qian, Chen and He, Ying},
	month = jan,
	year = {2024},
}

@misc{mo_efficient_2024,
	title = {Efficient {3D} {Shape} {Generation} via {Diffusion} {Mamba} with {Bidirectional} {SSMs}},
	url = {https://arxiv.org/abs/2406.05038v1},
	language = {en},
	urldate = {2026-06-03},
	journal = {arXiv.org},
	author = {Mo, Shentong},
	month = jun,
	year = {2024},
}

@misc{caytuiro_3d_2025,
	title = {{3D} {Shape} {Generation}: {A} {Survey}},
	shorttitle = {{3D} {Shape} {Generation}},
	url = {https://arxiv.org/abs/2506.22678v2},
	language = {en},
	urldate = {2026-06-03},
	journal = {arXiv.org},
	author = {Caytuiro, Nicolas and Sipiran, Ivan},
	month = jun,
	year = {2025},
}

@misc{jun_shap-e_2023,
	title = {Shap-{E}: {Generating} {Conditional} {3D} {Implicit} {Functions}},
	shorttitle = {Shap-{E}},
	url = {https://arxiv.org/abs/2305.02463v1},
	language = {en},
	urldate = {2026-06-03},
	journal = {arXiv.org},
	author = {Jun, Heewoo and Nichol, Alex},
	month = may,
	year = {2023},
}

@article{regenwetter_deep_2022,
	title = {Deep {Generative} {Models} in {Engineering} {Design}: {A} {Review}},
	volume = {144},
	issn = {1050-0472, 1528-9001},
	shorttitle = {Deep {Generative} {Models} in {Engineering} {Design}},
	url = {https://asmedigitalcollection.asme.org/mechanicaldesign/article/144/7/071704/1136676/Deep-Generative-Models-in-Engineering-Design-A},
	doi = {10.1115/1.4053859},
	language = {en},
	number = {7},
	urldate = {2026-06-03},
	journal = {Journal of Mechanical Design},
	author = {Regenwetter, Lyle and Nobari, Amin Heyrani and Ahmed, Faez},
	month = jul,
	year = {2022},
	pages = {071704},
}

@misc{liu_meshdiffusion_2023,
	title = {{MeshDiffusion}: {Score}-based {Generative} {3D} {Mesh} {Modeling}},
	shorttitle = {{MeshDiffusion}},
	url = {https://arxiv.org/abs/2303.08133v2},
	language = {en},
	urldate = {2026-06-03},
	journal = {arXiv.org},
	author = {Liu, Zhen and Feng, Yao and Black, Michael J. and Nowrouzezahrai, Derek and Paull, Liam and Liu, Weiyang},
	month = mar,
	year = {2023},
}

@misc{shi_deep_2023,
	title = {Deep {Generative} {Models} on {3D} {Representations}: {A} {Survey}},
	shorttitle = {Deep {Generative} {Models} on {3D} {Representations}},
	url = {http://arxiv.org/abs/2210.15663},
	doi = {10.48550/arXiv.2210.15663},
	urldate = {2026-06-03},
	publisher = {arXiv},
	author = {Shi, Zifan and Peng, Sida and Xu, Yinghao and Geiger, Andreas and Liao, Yiyi and Shen, Yujun},
	month = aug,
	year = {2023},
	note = {arXiv:2210.15663},
}

@inproceedings{siddiqui_meshgpt_2024,
	title = {{MeshGPT}: {Generating} {Triangle} {Meshes} with {Decoder}-{Only} {Transformers}},
	shorttitle = {{MeshGPT}},
	url = {https://openaccess.thecvf.com/content/CVPR2024/html/Siddiqui_MeshGPT_Generating_Triangle_Meshes_with_Decoder-Only_Transformers_CVPR_2024_paper.html},
        booktitle = {Proceedings of the IEEE/CVF Conference on Computer Vision and Pattern Recognition (CVPR)},
	language = {en},
	urldate = {2026-06-03},
	author = {Siddiqui, Yawar and Alliegro, Antonio and Artemov, Alexey and Tommasi, Tatiana and Sirigatti, Daniele and Rosov, Vladislav and Dai, Angela and Nießner, Matthias},
	year = {2024},
	pages = {19615--19625},
}

@misc{shi_mesh-informed_2025,
	title = {Mesh-{Informed} {Neural} {Operator} : {A} {Transformer} {Generative} {Approach}},
	shorttitle = {Mesh-{Informed} {Neural} {Operator}},
	url = {https://arxiv.org/abs/2506.16656v3},
	language = {en},
	urldate = {2026-06-03},
	journal = {arXiv.org},
	author = {Shi, Yaozhong and Ross, Zachary E. and Asimaki, Domniki and Azizzadenesheli, Kamyar},
	month = jun,
	year = {2025},
}

@article{weller_tensorial_1998,
	title = {A tensorial approach to computational continuum mechanics using object-oriented techniques},
	volume = {12},
	issn = {0894-1866},
	url = {https://pubs.aip.org/cip/article/12/6/620/510187/A-tensorial-approach-to-computational-continuum},
	doi = {10.1063/1.168744},
	language = {en},
	number = {6},
	urldate = {2026-01-31},
	journal = {Computers in Physics},
	author = {Weller, H. G. and Tabor, G. and Jasak, H. and Fureby, C.},
	month = nov,
	year = {1998},
	pages = {620--631},
}

@inproceedings{turk_zippered_1994,
	address = {Not Known},
	title = {Zippered polygon meshes from range images},
	copyright = {https://www.acm.org/publications/policies/copyright\_policy\#Background},
	isbn = {9780897916677},
	url = {http://portal.acm.org/citation.cfm?doid=192161.192241},
	doi = {10.1145/192161.192241},
	language = {en},
	urldate = {2026-01-31},
	booktitle = {Proceedings of the 21st annual conference on {Computer} graphics and interactive techniques  - {SIGGRAPH} '94},
	publisher = {ACM Press},
	author = {Turk, Greg and Levoy, Marc},
	year = {1994},
	pages = {311--318},
}

@inproceedings{potamias_shapefusion_2025,
	address = {Cham},
	title = {{ShapeFusion}: {A} {3D} {Diffusion} {Model} for {Localized} {Shape} {Editing}},
	isbn = {9783031726309},
	shorttitle = {{ShapeFusion}},
	doi = {10.1007/978-3-031-72630-9_5},
	language = {en},
	booktitle = {Computer {Vision} – {ECCV} 2024},
	publisher = {Springer Nature Switzerland},
	author = {Potamias, Rolandos Alexandros and Tarasiou, Michail and Ploumpis, Stylianos and Zafeiriou, Stefanos},
	editor = {Leonardis, Aleš and Ricci, Elisa and Roth, Stefan and Russakovsky, Olga and Sattler, Torsten and Varol, Gül},
	year = {2025},
	pages = {72--89},
}

@article{xiong_octfusion_2025,
	title = {{OctFusion}: {Octree}‐based {Diffusion} {Models} for {3D} {Shape} {Generation}},
	volume = {44},
	issn = {0167-7055, 1467-8659},
	shorttitle = {{OctFusion}},
	url = {https://onlinelibrary.wiley.com/doi/10.1111/cgf.70198},
	doi = {10.1111/cgf.70198},
	language = {en},
	number = {5},
	urldate = {2026-01-30},
	journal = {Computer Graphics Forum},
	author = {Xiong, Bojun and Wei, Si‐Tong and Zheng, Xin‐Yang and Cao, Yan‐Pei and Lian, Zhouhui and Wang, Peng‐Shuai},
	month = aug,
	year = {2025},
	pages = {e70198},
}

@misc{padula_generative_2025,
	title = {Generative {Models} for {Parameter} {Space} {Reduction} applied to {Reduced} {Order} {Modelling}},
	url = {http://arxiv.org/abs/2506.09721},
	doi = {10.48550/arXiv.2506.09721},
	urldate = {2026-01-30},
	publisher = {arXiv},
	author = {Padula, Guglielmo and Rozza, Gianluigi},
	month = jun,
	year = {2025},
	note = {arXiv:2506.09721},
}

@article{sederberg_free-form_1986,
	title = {Free-form deformation of solid geometric models},
	volume = {20},
	doi = {10.1145/15886.15903},
	number = {4},
	journal = {ACM SIGGRAPH Computer Graphics},
	author = {Sederberg, Thomas W. and Parry, Scott R.},
	year = {1986},
	pages = {151--160},
}

@article{chapelier_free-form_2021,
	title = {Free-{Form} {Deformation} {Digital} {Image} {Correlation} ({FFD}-{DIC}): {A} non-invasive spline regularization for arbitrary finite element measurements},
	volume = {384},
	doi = {10.1016/j.cma.2021.113992},
	journal = {Computer Methods in Applied Mechanics and Engineering},
	author = {Chapelier, M. and Bouclier, R. and Passieux, J.-C.},
	year = {2021},
	pages = {113992},
}

@article{padula_generative_2024,
	title = {Generative models for the deformation of industrial shapes with linear geometric constraints: {Model} order and parameter space reductions},
	volume = {423},
	doi = {10.1016/j.cma.2024.116823},
	journal = {Computer Methods in Applied Mechanics and Engineering},
	author = {Padula, Guglielmo and Romor, Francesco and Stabile, Giovanni and Rozza, Gianluigi},
	year = {2024},
	pages = {116823},
}

@article{fukuda_efficient_2024,
	title = {Efficient musculoskeletal annotation using free-form deformation},
	volume = {14},
	doi = {10.1038/s41598-024-67125-3},
	number = {1},
	journal = {Scientific Reports},
	author = {Fukuda, Norio and Konda, Shoji and Umehara, Jun and Hirashima, Masaya},
	year = {2024},
	pages = {16077},
}

@article{wang_diffusion_2025,
	title = {Diffusion {Models} for {3D} {Generation}: {A} {Survey}},
	volume = {11},
	doi = {10.26599/CVM.2025.9450452},
	number = {1},
	journal = {Computational Visual Media},
	author = {Wang, Chen and Peng, Hao-Yang and Liu, Ying-Tian and Gu, Jiatao and Hu, Shi-Min},
	year = {2025},
	pages = {1--28},
}

@article{liu_novel_2025,
	title = {A {Novel} body-fitted free form deformation method for hull form optimization},
	volume = {342},
	doi = {10.1016/j.oceaneng.2025.122786},
	journal = {Ocean Engineering},
	author = {Liu, Xinwang and Rong, Zitong and Yuan, Lei and Ji, Xiaohang and Wang, Luyao and Sun, Xu},
	year = {2025},
	pages = {122786},
}

@article{scroggs_basix_2022,
	title = {Basix: a runtime finite element basis evaluation library},
	volume = {7},
	issn = {2475-9066},
	doi = {10.21105/joss.03982},
	number = {73},
	journal = {Journal of Open Source Software},
	author = {Scroggs, Matthew W. and Baratta, Igor A. and Richardson, Chris N. and Wells, Garth N.},
	month = may,
	year = {2022},
	pages = {3982},
}

@article{baratta_dolfinx_2023,
	title = {{DOLFINx}: {The} next generation {FEniCS} problem solving environment},
	issn = {0098-3500},
	doi = {10.1145/3524456},
	journal = {ACM Iransactions on Mathematical Software (IUMS)}, 
	author = {Baratta, Igor A. and Dean, Joseph P. and Dokken, Jørgen S. and Habera, Michal and Hale, Jack S. and Richardson, Chris N. and Rognes, Marie E. and Scroggs, Matthew W. and Sime, Nathan and Wells, Garth N.},
	year = {2023},
    volume= {48},
    issue ={2},
}

@article{salmoiraghi_free-form_2018,
	title = {Free-form deformation, mesh morphing and reduced-order methods: enablers for efficient aerodynamic shape optimisation},
	volume = {32},
	issn = {1061-8562},
	doi = {10.1080/10618562.2018.1514115},
	number = {4-5},
	journal = {International Journal of Computational Fluid Dynamics},
	author = {Salmoiraghi, F. and Scardigli, A. and Telib, H. and Rozza, G.},
	month = may,
	year = {2018},
	pages = {233--247},
}

@article{grey_active_2018,
	title = {Active {Subspaces} of {Airfoil} {Shape} {Parameterizations}},
	volume = {56},
	issn = {0001-1452},
	doi = {10.2514/1.J056054},
	number = {5},
	journal = {AIAA Journal},
	author = {Grey, Zachary J. and Constantine, Paul G.},
	month = may,
	year = {2018},
	pages = {2003--2017},
}

@article{padula_brief_2024,
	title = {A brief review of reduced order models using intrusive and non‐intrusive techniques},
	volume = {24},
	issn = {1617-7061},
	doi = {10.1002/pamm.202400210},
	number = {4},
	journal = {PAMM},
	author = {Padula, Guglielmo and Girfoglio, Michele and Rozza, Gianlugi},
	month = dec,
	year = {2024},
}

@book{constantine_active_2015,
	address = {Philadelphia, PA},
	title = {Active {Subspaces}: {Emerging} {Ideas} for {Dimension} {Reduction} in {Parameter} {Studies}},
	isbn = {9781611973853 9781611973860},
	shorttitle = {Active {Subspaces}},
	url = {http://epubs.siam.org/doi/book/10.1137/1.9781611973860},
	language = {en},
	urldate = {2026-06-16},
	publisher = {Society for Industrial and Applied Mathematics},
	author = {Constantine, Paul G.},
	month = mar,
	year = {2015},
	doi = {10.1137/1.9781611973860},
}

@article{romor_kernelbased_2022,
	title = {Kernel‐based active subspaces with application to computational fluid dynamics parametric problems using the discontinuous {Galerkin} method},
	volume = {123},
	issn = {0029-5981, 1097-0207},
	url = {https://onlinelibrary.wiley.com/doi/10.1002/nme.7099},
	doi = {10.1002/nme.7099},
	language = {en},
	number = {23},
	urldate = {2026-06-16},
	journal = {International Journal for Numerical Methods in Engineering},
	author = {Romor, Francesco and Tezzele, Marco and Lario, Andrea and Rozza, Gianluigi},
	month = dec,
	year = {2022},
	pages = {6000--6027},
}

@article{feng_parametric_2021,
	title = {Parametric Hull Form Optimization of Containerships for Minimum Resistance in Calm Water and in Waves},
	volume = {20},
	issn = {1993-5048},
	url = {https://doi.org/10.1007/s11804-021-00243-w},
	doi = {10.1007/s11804-021-00243-w},
	pages = {670--693},
	number = {4},
	journal = {Journal of Marine Science and Application},
	shortjournal = {J. Marine. Sci. Appl.},
    year={2021},
	author = {Feng, Yanxin and el Moctar, Ould and Schellin, Thomas E.},
	urldate = {2026-09-20},
	date = {2021-12-01},
	langid = {english},
}

@article{moctar_duisburg_2012,
	title = {Duisburg Test Case: Post-Panamax Container Ship for Benchmarking},
	volume = {59},
	issn = {0937-7255, 2056-7111},
	url = {http://www.tandfonline.com/doi/full/10.1179/str.2012.59.3.004},
	doi = {10.1179/str.2012.59.3.004},
	shorttitle = {Duisburg Test Case},
	pages = {50--64},
    year = {2012},
	number = {3},
	journal = {Ship Technology Research},
	shortjournal = {Ship Technology Research},
	author = {Moctar, Ould El and Shigunov, Vladimir and Zorn, Tobias},
	urldate = {2026-09-20},
	date = {2012-08},
	langid = {english},
}
\appendix

\section{Numerical set-up of the DTC hull test case}
\label{app:cfd}

The hydrodynamic simulations of Section~\ref{sec:testcases} are based on the
\texttt{DTCHull} tutorial distributed with OpenFOAM~v2412, modified only in the bow
region through the free-form deformation described in Section~\ref{sec:testcases}. For
reproducibility, and since the deformation alters the geometry on which the drag is
evaluated, the complete set-up is summarised here.

\subsection{Geometry and scale}
\label{app:cfd:geometry}

The Duisburg Test Case is a post-panamax container ship
benchmark~\cite{moctar_duisburg_2012}. The simulations are performed at model scale
$1{:}59.407$, at which the principal dimensions of the hull are those of
Table~\ref{tab:dtc_dimensions}. The still-water plane is located at $z=0.244$~m above
the keel, so that the undeformed hull floats at its design draught. Only half of the
hull and of the surrounding domain is discretised, the centreplane $y=0$ being treated
as a plane of symmetry; all drag values reported in Section~\ref{sec:testcases}
therefore refer to half the hull and must be doubled for comparison with full-hull
measurements.

\begin{table}[ht]
\caption{Principal dimensions of the DTC hull, at full scale and at the model scale
$1{:}59.407$ used here~\cite{moctar_duisburg_2012, feng_parametric_2021}.}
\label{tab:dtc_dimensions}
\centering
\begin{tabular}{llcc}
\hline
Quantity & Symbol & Full scale & Model scale \\
\hline
Length between perpendiculars & $L_{pp}$ [m]      & $355.0$    & $5.976$ \\
Beam                          & $B$ [m]           & $51.0$     & $0.859$ \\
Design draught                & $T$ [m]           & $14.5$     & $0.244$ \\
Displacement                  & $\nabla$ [m$^{3}$]& $173\,467$ & $0.827$ \\
Wetted surface area           & $S$ [m$^{2}$]     & $22\,032$  & $6.243$ \\
Block coefficient             & $C_{B}$           & \multicolumn{2}{c}{$0.661$} \\
Design speed                  & [kn]              & $25.0$     & $3.244$ \\
\hline
\end{tabular}
\end{table}

The free-form deformation acts on the bulbous bow, within a lattice box of origin
$(5.80,\,-0.10,\,-0.05)$~m and edge lengths $(0.50,\,0.10,\,0.33)$~m, that is on a
region of streamwise extent $0.084\,L_{pp}$ located at the forward perpendicular. The
remainder of the hull is left unmodified, and the control-point weights are projected
onto the constant-volume manifold as described in Section~\ref{sec:testcases}, so that
every deformed configuration encloses the same volume as the reference geometry to
machine precision. The displacement of Table~\ref{tab:dtc_dimensions} is therefore
common to the whole dataset, and differences in drag between configurations cannot be
attributed to changes in immersed volume.

\subsection{Flow conditions and fluid properties}
\label{app:cfd:flow}

The hull is held fixed and the flow is imposed at the inlet with uniform speed
$U_{0}=1.668$~m\,s$^{-1}$ directed along the negative $x$-axis. This is the design
speed of the benchmark, $25$~kn at full scale, and corresponds to
\[
Fn=\frac{U_{0}}{\sqrt{g L_{pp}}}=0.218 ,
\qquad
Re=\frac{U_{0}L_{pp}}{\nu_{W}}=9.1\times10^{6},
\]
with $g=9.81$~m\,s$^{-2}$. The two phases are treated as isothermal, immiscible and
Newtonian, with the properties listed in Table~\ref{tab:fluid_properties}; surface
tension is neglected, which is standard practice for ship-resistance computations at
this scale, where the Weber number is large.

\begin{table}[ht]
\caption{Fluid properties. Kinematic viscosities and densities are those of the
OpenFOAM \texttt{DTCHull} tutorial.}
\label{tab:fluid_properties}
\centering
\begin{tabular}{lcc}
\hline
Quantity & Water & Air \\
\hline
Density $\rho$ [kg\,m$^{-3}$]              & $998.8$    & $1.0$ \\
Kinematic viscosity $\nu$ [m$^{2}$\,s$^{-1}$] & $1.09\times10^{-6}$ & $1.48\times10^{-5}$ \\
\hline
Surface tension $\sigma$ [N\,m$^{-1}$] & \multicolumn{2}{c}{$0$} \\
\hline
\end{tabular}
\end{table}

\subsection{Computational domain and boundary conditions}
\label{app:cfd:bcs}

The domain is a rectangular box extending $7.0\,L_{pp}$ in the streamwise direction,
$3.2\,L_{pp}$ in the transverse direction from the centreplane, and $2.7\,L_{pp}$ below
the still-water plane, with a further $3.76$~m of air above it. The inlet lies
$1.7\,L_{pp}$ ahead of the bow and the outlet $4.4\,L_{pp}$ astern, distances
sufficient to avoid reflection of the generated wave system at the boundaries within
the integration horizon. The boundary conditions are listed in
Table~\ref{tab:bcs}. The volume fraction is initialised to $\alpha=1$ below $z=0.244$~m
and $\alpha=0$ above it, both in the interior and on the boundary faces, so that the
inlet delivers water below the waterline and air above it.

\begin{table}[ht]
\caption{Boundary conditions. The hull is stationary, so the \texttt{movingWallVelocity}
condition reduces to a no-slip condition.}
\label{tab:bcs}
\centering
\small
\renewcommand{\arraystretch}{1.15}
\begin{tabularx}{\linewidth}{>{\raggedright\arraybackslash}p{2.2cm}
                              >{\raggedright\arraybackslash}X
                              >{\raggedright\arraybackslash}X
                              >{\raggedright\arraybackslash}X
                              >{\raggedright\arraybackslash}X}
\hline
Patch & $\mathbf{u}$ & $p_{rgh}$ & $\alpha$ & Turbulence \\
\hline
Inlet      & fixed value $(-U_{0},0,0)$ & fixed flux    & fixed value    & fixed value \\
Outlet     & outlet phase mean velocity & zero gradient & variable height flow rate & inlet--outlet \\
Atmosphere & pressure inlet--outlet     & total pressure $p_{0}=0$ & inlet--outlet & inlet--outlet \\
Hull       & no slip                    & fixed flux    & zero gradient  & wall functions \\
Centreplane, bottom, side & \multicolumn{4}{l}{symmetry plane} \\
\hline
\end{tabularx}
\end{table}
At the inlet the turbulence quantities are set to $k=1.5\times10^{-4}$~m$^{2}$\,s$^{-2}$
and $\omega=2$~s$^{-1}$, corresponding to a turbulence intensity of approximately
$0.6\%$ and an eddy viscosity ratio of order unity.

\subsection{Discretisation and mesh quality}
\label{app:cfd:mesh}

The background grid is generated with \texttt{blockMesh} and comprises $134\,064$
hexahedral cells, vertically graded so that the cell height in the neighbourhood of the
still-water plane is approximately $1.5\times10^{-2}$~m, that is $T/16$. Six successive
applications of \texttt{refineMesh} then refine the grid in the streamwise and
transverse directions only, within nested boxes of decreasing size centred on the hull
and on the free surface, the innermost of which spans
$(-0.5,-0.55,-0.15)$--$(6.25,0,0.65)$~m. This brings the horizontal cell size in the
region of interest to approximately $1.6\times10^{-2}$~m, comparable to the vertical
one, so that the cells resolving the free surface and the hull are nearly isotropic and
the interface is captured across roughly one cell. The hull surface is then introduced
by \texttt{snappyHexMesh}, with three near-wall layers of expansion ratio $1.5$ and
final-layer thickness equal to $0.7$ of the adjacent cell size. The grid of the
reference configuration comprises $845\,539$ cells and $919\,139$ points. The mesh is
regenerated for every deformed configuration, the deformation being applied to the
surface triangulation before \texttt{snappyHexMesh} is run, so that the near-wall
resolution is preserved across the dataset; since the snapping is performed
independently for each geometry, the cell and point counts fluctuate mildly about these
reference values.

Geometries produced by the generative model are checked for validity before meshing:
the decoded weights are subjected to the same volume projection as the sampled ones,
and the resulting surfaces are verified to be admissible inputs to the meshing stage.
Mesh quality is then monitored across the whole dataset with \texttt{checkMesh}. For
every configuration the maximum face non-orthogonality is $70^{\circ}$ and the mean is
$7^{\circ}$, and at most a single face out of approximately $8.5\times10^{5}$ cells is
flagged as highly skew, with a maximum skewness of $6$. These figures are the same
across the deformed configurations and as those of the undeformed reference mesh, so
the deformation does not degrade the discretisation: the mean non-orthogonality, which
governs the accuracy of the pressure equation over the bulk of the domain, is an order
of magnitude below the tolerance of the utility, and an isolated skew face contributes
a negligible share of the hull-surface integral. Differences in drag between
configurations are therefore not attributable to cells created by the deformation.
\begin{figure}
    \centering
\includegraphics[width=0.9\linewidth]{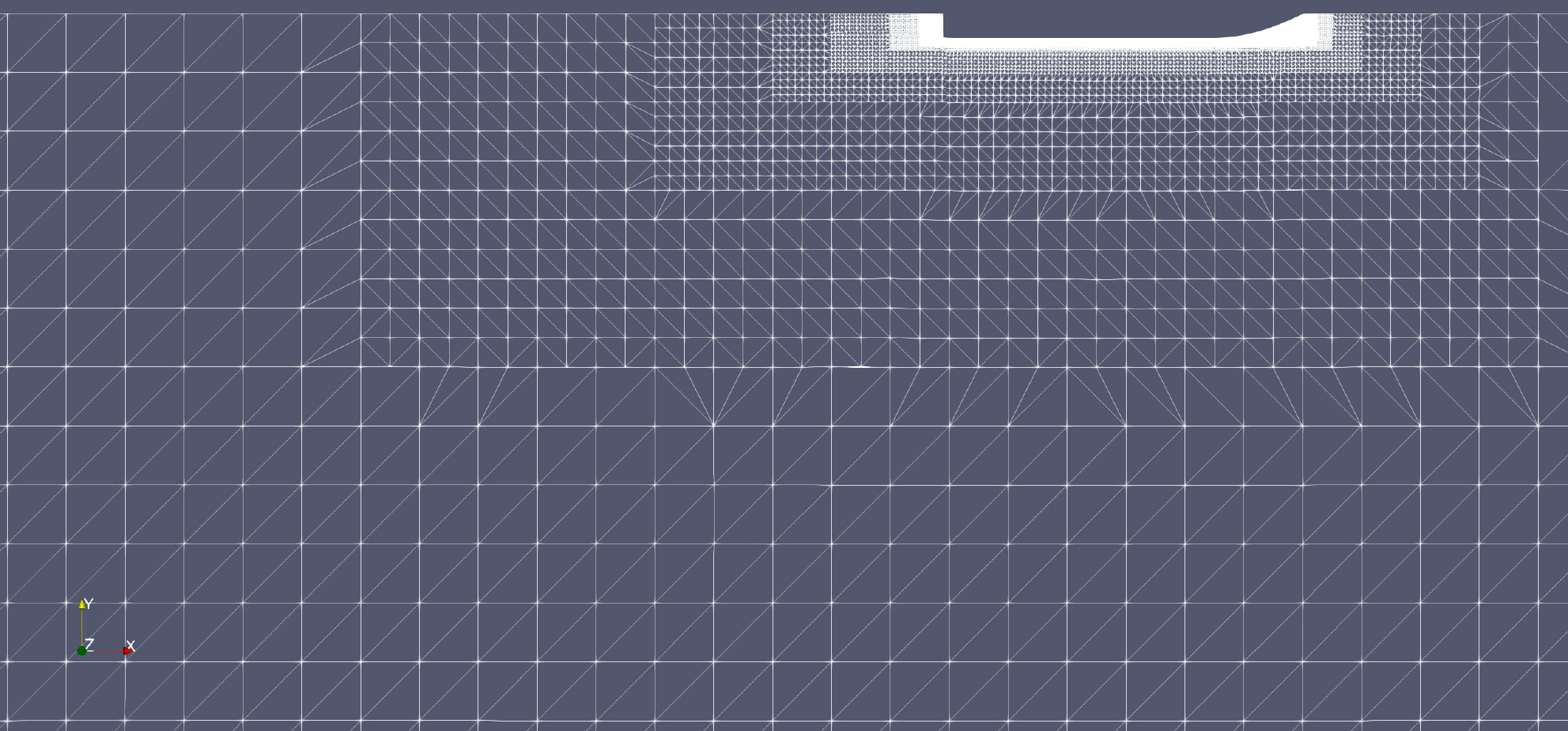}
    \caption{A top view of the computational mesh with the six \texttt{refineMesh} functions applied.}
\end{figure}

\subsection{Turbulence model and wall treatment}
\label{app:cfd:turbulence}

Turbulence is modelled with the $k$--$\omega$ Shear Stress Transport model in its
standard OpenFOAM implementation. The hull is treated as a rough wall through the
\texttt{nutkRoughWallFunction} condition, with equivalent sand-grain roughness
$K_{s}=100\,\mu$m and roughness constant $C_{s}=0.5$, together with wall functions for
$k$ and $\omega$; the values quoted below are accordingly those of the rough-wall
formulation. The boundary layer is not resolved down to the viscous sublayer. The
dimensionless wall distance computed from the converged solution has a mean of $70$
over the hull, with extreme values of $3$ and $1016$. The mean lies well within the
range of validity of the log-law formulation, as required by the wall-function
treatment. The extremes are localised: the lowest values occur near stagnation and
separation points, where the friction velocity vanishes and $y^{+}$ is necessarily
small, and the highest on the part of the hull lying above the waterline, where the
cells are filled with air and the reported value, computed with the properties of the
lighter phase, has no bearing on the resistance.

\subsection{Time integration and convergence}
\label{app:cfd:time}

Steady solutions are sought with the \texttt{interFoam} solver used in local
time-stepping mode: the temporal derivative is discretised with a local Euler scheme,
so that each cell advances at its own rate and the transient has no physical meaning,
the iteration index acting as a pseudo-time. The local step is bounded by a maximum
Courant number of $10$, a maximum interface Courant number of $5$ and a maximum step of
unity, with the usual smoothing and damping of the reciprocal step field. Pressure and
velocity are coupled with the PIMPLE algorithm using two pressure correctors and no
momentum predictor, while the volume fraction is advanced with a semi-implicit MULES
limiter and two corrector steps.

Each simulation is advanced for $4000$ pseudo-time iterations, at which point the
integrated force on the hull has settled. The drag is taken as the streamwise component
of the pressure and viscous forces on the hull patch, evaluated at every iteration by
the \texttt{forces} function object, and is averaged over the last $500$ iterations to
remove the residual oscillation associated with the free-surface treatment; over that
window the force fluctuates by approximately $1\%$ of its mean. The comparison of the
resulting baseline resistance with towing-tank and published numerical data is reported
in Section~\ref{sec:testcases}.

\subsection{Software}
\label{app:cfd:software}

All simulations are performed with OpenFOAM~v2412~\cite{weller_tensorial_1998}, using 8 processors. Mesh
generation uses the \texttt{blockMesh}, \texttt{surfaceFeatures}, \texttt{topoSet},
\texttt{refineMesh} and \texttt{snappyHexMesh} utilities of the same release, and the
free-form deformation is applied to the surface triangulation before meshing.

\end{document}